# Beyond Gödel

Bhupinder Singh Anand

## Simply consistent constructive systems of first order Peano's Arithmetic that do not yield undecidable propositions by Gödel's reasoning

**In this paper, we consider significant alternative systems of first order predicate calculus. We consider systems where the meta-assertion "PA proves: $F(x)$" translates under interpretation as "$F(x)$ is satisfied for all values of $x$, in the domain of the interpretation, that can be formally represented in PA", whilst "PA proves: $(Ax)F(x)$" interprets as "$F(x)$ is satisfied by all values of $x$ in the domain of the interpretation".**

**We are thus able to argue that formal systems of first order Arithmetic that admit Gödelian undecidable propositions validly are abnormally non-constructive.**

**We argue that, in such systems, the strong representation of primitive recursive predicates admits abnormally non-constructive, Platonistic, elements into the formal system that are not reflected in the predicates which they are intended to formalise.**

**We argue that the source of such abnormal Platonistic elements in these systems is the non-constructive Generalisation rule of inference of first order logic.**

**We argue that, in most simply consistent systems that faithfully formalise intuitive Arithmetic, we cannot infer from Gödel's reasoning the Platonistic existence of abnormally non-constructive propositions that are formally undecidable, but true under interpretation.**

We define a constructive formal system of Peano's Arithmetic, *omega*$_2$-**PA**, whose axioms are identical to the axioms of standard Peano's Arithmetic **PA**, but lead to significantly different logical consequences.

We thus argue that the formal undecidability of true Arithmetical propositions is a characteristic not of relations that are Platonistically inherent in any Arithmetic of the natural numbers, but of the particular formalisation chosen to represent them.

We finally argue that Gödel's reasoning essentially formalises the *Liar* sentence in standard **PA** by means of a "palimpsest".

## Chapter 1. First order Arithmetic is not *omega*-constructive

**Contents**

## 1.0. Introduction

We take as our primary reference Mendelson's first order exposition of the essentially second order formal system and various revolutionary concepts introduced by Gödel in his original paper "*On formally undecidable propositions of Principia Mathematica and related systems I*".

We also borrow some content, and style of presentation, from *Karlis Podnieks*' proof of *Goedel's incompleteness theorem* in his e-textbook "*Around Goedel's Theorem*".[1]

We take as our starting point the systems of first order predicate calculus and first order Arithmetic defined by Mendelson (*p57 & p102*).

We introduce the concept of *omega*-constructivity, and show that standard first order Arithmetic **PA** is not *omega*-constructive.

---

[1] Where possible, words and phrases have been hyper-linked to references, with detailed locations or further hyperlinks to the Internet, for consistency in definitions and meanings, or for expanding on their possible meanings-in-usage in a wider context.

Pointing your cursor at a word or phrase will indicate if the word or phrase has been hyper-linked to another reference (the cursor should change shape in most browsers).

Clicking once (default used whilst composing this paper) will bring up the cross-linked reference on the screen.

Use the 'Back' button '<=' on the browser toolbar at the top of the screen (not any of the '**?** ' buttons in the document) to return to the original hyper-linked word or phrase in the text.

### 1.1. Introducing *omega*-constructivity : *omega*-PA

We say that a simply consistent formal system of standard first order Arithmetic **PA** is *omega*-constructive if and only if it is not the case that, for any well-formed formula $F(x_1, x_2, x_3, \ldots, x_n)$ of **PA**, we can have both:

(*i*)   ~(**PA** proves: $F(x_1, x_2, x_3, \ldots, x_n)$), and

(*ii*)  $(Ax_1, Ax_2, Ax_3, \ldots, Ax_n)$(**PA** proves: $F(x_1, x_2, x_3, \ldots, x_n)$), where the domain of $x_i$ is the set of natural numbers.

In other words, we cannot have a situation where a well-formed formula such as $F(x)$ is unprovable in **PA**, but $F(1), F(2), F(3), \ldots$ are all provable, if **PA** is both simply consistent and *omega*-constructive, where $n$ denotes the numeral representing[2] the natural number $n$ in the formal system **PA**.

We name such an *omega*-constructive **PA** system 'first order *omega*-Arithmetic', and refer to it as *omega*-**PA**.

### 1.2. The significance of *omega*-constructivity

This definition implicitly reflects the thesis that an infinite, compound, meta-assertion such as "$F(1)$ is provable in **PA** & $F(2)$ is provable in **PA** & $F(3)$ is provable in **PA** & ..." is simply a convenient way of expressing and visualising the meta-assertion "$F(n)$ is provable in **PA** for all $n$", where $n$ is intuitively taken to range over the natural numbers.

---

[2] See also Gödel (*Gödel 1931 Def17*).

### 1.3. Constructive PA should be intuitive

We argue that, in a constructive formal system of **PA** that is intended to formalise the natural numbers - where the intended domain of $x$ under interpretation contains, and preferably consists only of, elements that correspond to the numerals that represent the natural numbers in the formal system - the meta-assertion "$F(x)$ is provable in **PA**" should intuitively translate as equivalent to "$F(n)$ is provable in **PA** for all $n$", where $n$ is intuitively taken to range over the natural numbers.

So a formal system that is not *omega*-constructive is counter-intuitive.

### 1.4. Differentiating between *omega*-constructivity and *omega*-completeness

We note the difference between *omega*-constructivity and *omega*-completeness.

**PA** is said to be *omega*-complete if and only if there is a no well-formed formula $F(x)$ with one free variable such that:

(i)   $\sim$( **PA** proves: $(Ax)F(x)$), and

(ii)  $(An)$( **PA** proves: $F(n)$), where the domain of $n$ is the set of natural numbers.

Now, in our intended formal system **PA** of the natural numbers, an infinite, compound, meta-assertion such as "$F(1)$ is provable in **PA** & $F(2)$ is provable in **PA** & $F(3)$ is provable in **PA** & ... " can also be intuitively asserted as a consequence of a meta-proof of the generalization "$(Ax)(F(x)$ is provable in **PA**)".

However, we need to consider the possibility that, under interpretation, such a **PA** may admit an infinite domain for $x$ that is not denumerable, or that the generalisation "$(Ax)(F(x)$" may be a well-formed formula, but an ill-defined proposition under interpretation (*which could be the case if **F**(x) either represents, or translates under interpretation as, a vague predicate*).

It follows that, under our intended formalisation, the denumerable meta-assertion "**F**(*1*) is provable in **PA** & **F**(*2*) is provable in **PA** & **F**(*3*) is provable in **PA** & ..." cannot be assumed to conversely yield a proof of '(A$x$)**F**($x$)' in **PA**.

So a system that is not *omega*-complete is not necessarily counter-intuitive.

**1.5. The relevance of *omega*-constructivity**

To see the relevance of *omega*-constructivity to **PA**, we consider in §**1.5** to §**1.7** the primitive recursive (*Mendelson p120*) predicate ***prf'***($u$, $y$), which we define[3] as true if and only if $u$ is the Gödel-number (*Podnieks 2001 Section 5.2*) of a well-formed formula **F**($x_1$, $x_2$, $x_3$, ... , $x_n$) of **PA**, containing the free variables $x_1$, $x_2$, $x_3$, ... , $x_n$, and $y$ is the Gödel-number of a proof in **PA** of **F**(***u***, $x_2$, $x_3$, ... , $x_n$).

**1.6. *prf'*($x$, $y$) is a Turing decidable predicate**

Now we can construct a Turing machine that, given any $u$ and $v$, will decode (*Podnieks 2001 Excercise 5.3*) $u$ into **F**($x_1$, $x_2$, $x_3$, ... , $x_n$), check whether $x_1$ occurs in it, construct **F**(***u***, $x_2$, $x_3$, ... , $x_n$) by replacing $x_1$ wherever it occurs by ***u***, decode $v$, check if $v$ is a proof sequence, and finally check to see if **F**(***u***, $x_2$, $x_3$, ... , $x_n$) is the last well-formed formula in the proof sequence coded by $v$.

Thus ***prf'***($x$, $y$) is a Turing-decidable predicate.

---

[3] This definition is based on the predicate $W_1$($u$,$y$) used by Mendelson (*p143*) in his proof of Gödel's undecidability theorem. It corresponds to Gödel's (*1931 p188 eqn 8.1*) primitive recursive predicate $x$B(Sb($y$ : (19, Z($y$)))).

## 1.7. The Representation Theorem

By Podnieks Representation Theorem (*Section 3.3*) we can therefore construct a **PA**-formula **PRF'**($x$, $y$) expressing the predicate **prf'**($x$, $y$) such that, for any given $k_1$, $k_2$:

(i)  **prf'**($k_1$, $k_2$) is true  => **PA** proves: **PRF'**($k_1$, $k_2$), and

(ii) ~**prf'**($k_1$, $k_2$) is true => **PA** proves: ~**PRF'**($k_1$, $k_2$).

We let $r$ be the Gödel-number of ~**PRF'**($x_1$, $x_2$).

We consider the well-formed formula ~**PRF'**($r$, $x_2$).

Now, from the definition of **prf'**($u$, $y$), it follows that[4]:

(iii) **prf'**($r$, $y$) is true  <=> $y$ is the Gödel-number of a proof in **PA**
of ~**PRF'**($r$, $x_2$).

## 1.8. Is PA both simply consistent and *omega*-constructive?

In §**1.8** to §**1.11**, we next address the question: Can **PA** be both simply consistent and *omega*-constructive?

Let us assume that it is, and that *omega*-**PA** proves ~**PRF'**($r$, $x_2$), so that some natural number $k$ is the Gödel-number of this proof.

Then, by §**1.7**(iii), **prf'**($r$, $k$) is true and hence,

   *omega*-**PA** proves: **PRF'**($r$, $k$).

However, by our premise, we have that,

---

[4] This meta-mathematical equivalence seems to be the counter-part of Podnieks formal self-reference lemma.

omega-**PA** proves: **~PRF'**(*r*, $x_2$), and so

omega-**PA** proves: **~PRF'**(*r*, *k*).

Therefore, since *omega*-**PA** is assumed simply consistent, it follows that **~PRF'**(*r*, $x_2$) cannot be proved in *omega*-**PA**.

### 1.9. Consequences of the *omega*-constructivity of *omega*-**PA**

We thus have from the *omega*-constructivity of *omega*-**PA** that it is not the case that:

(A*n*)(*omega*-**PA** proves: **~PRF'**(*r*, *n*)).

Now, since ***prf'***(*x*, *y*) is primitive recursive, if **~*prf'***(*r*, *y*) is true for all *y*, it follows meta-mathematically from §**1.7**(ii) that:

(A*n*)(*omega*-**PA** proves: **~PRF'**(*r*, *n*)).

So, if it is not the case that,

(A*n*)(*omega*-**PA** proves: **~PRF'**(*r*, *n*)),

then it is not the case that,

**~*prf'***(*r*, *y*) is true for all *y*.

Hence we have, meta-mathematically, that:

**~*prf'***(*r*, *k*) is false for some *k*,

***prf'***(*r*, *k*) is true,

*omega*-**PA** proves: **~PRF'**(*r*, $x_2$).

## 1.10. The meta-mathematical inconsistency in *omega*-PA

We have thus demonstrated that if **PA** is assumed both simply consistent and *omega*-constructive, then we have the meta-mathematical inconsistency:

    (*i*)     *omega*-**PA** proves: *~PRF'*(*r*, $x_2$)

                      =>   ~(*omega*-**PA** proves: *~PRF'*(*r*, $x_2$)), and

    (*ii*)   ~(*omega*-**PA** proves: *~PRF'*(*r*, $x_2$))

                      =>   *omega*-**PA** proves: *~PRF'*(*r*, $x_2$).

## 1.11. A simply consistent PA is neither *omega*-constructive, nor *omega*-complete

So, if standard first order Arithmetic **PA** is assumed simply consistent, then it cannot be *omega*-constructive in the sense defined above in §**1.1**.

Since, in standard first order Arithmetic **PA**, we have that:

    **PA** proves: *F*(*x*)    <=>    **PA** proves: (A*x*)*F*(*x*),

it follows that standard first order Arithmetic **PA** is, additionally, not *omega*-complete.

## 1.12. What is the significance of the *omega*-incompleteness of PA?

Now, if standard first order Arithmetic **PA** is assumed simply consistent, then there is the well-formed formula *~PRF'*(*r*, $x_2$) with one free variable such that:

    (*i*)    ~(**PA** proves: (A$x_2$)*~PRF'*(*r*, $x_2$)), and

    (*ii*)   (A*n*)(**PA** proves: *~PRF'*(*r*, *n*)).

The question arises: What is the significance of such *omega*-incompleteness of **PA**?

To get some indirect perspective on this question, let us for a moment consider the Gödelian argument.

### 1.13. The Gödelian argument

One interpretation, of the Gödelian argument as offered by J.R.Lucas, is the thesis that there is some non-mechanistic element - knowledge of which is Platonistically available to human intelligence but cannot be reflected in any machine intelligence - that is involved in establishing that a well-formed formula such as $(x)F(x)$ is formally unprovable in **PA**, yet translates into a true proposition under interpretation in the standard model.

However, if the foregoing arguments are simply non-constructive and intuitionistically objectionable, but not non-mechanistic, then the Platonistic elements of the reasoning are built into the formal logic itself, and can be built equally rigorously, independent of any model, into machine intelligence.

### 1.14. *omega*-provable propositions in PA

Thus if we introduce:

    **~**(**PA** proves: $(Ax)F(x)$)  &  $(An)$(**PA** proves: $F(n)$)

as the meta-definition of:

    '$(Ax)F(x)$' is *omega*-provable in **PA**,

we can reasonably argue that **PA** establishes the meta-assertions:

    '$(Ax_2)$**~*PRF*'**$(r, x_2)$' is *omega*-provable in **PA**,

    '**~*PRF*'**$(r, n)$' is provable in **PA** for every natural number $n$,

'**~PRF′**(*r*, *n*)' is true for every natural number *n* in every model of **PA**,

'(A*n*)**~PRF′**(*r*, *n*)' is true in every model of **PA**,

although **PA** is unable to prove (A*x*$_2$)**~PRF′**(***r***, *x*$_2$), and even though (A*x*)**~PRF′**(*r*, *x*) may not hold true in every model of **PA**.

### 1.15.  What kind of propositions are *omega*-provable?

The question arises: What makes a proposition such as (A*x*)***F***(*x*) *omega*-provable?

An obvious possibility, as remarked earlier, is that ***F***(*x*) may be a well-formed formula, but translate under interpretation into a vague or ill-defined predicate similar to the paradoxical self-referential constructions that give rise to the various linguistic, logical and mathematical antinomies.

We note that Gödel's aim in his 1931 paper "*On formally undecidable propositions of Principia Mathematica and related systems I*" was not to eliminate the antinomies, but to eliminate imprecision of expression in the concepts of truth and proof by rigorously formalising (*Gödel 1931 p176*) first order logic and Arithmetic.

Accordingly, in his 1934 lecture notes (*Davis p46*) "*On undecidable propositions of formal mathematical systems*", Gödel notes that "We shall depend on the theory of types as our means for avoiding paradox".

So, although ***F***(*x*) may be provable for any given *x*, the provability of (A*x*)***F***(*x*) may, in some interpretation, involve an implicit reference to an ill-defined totality of all *x* that satisfy ***F***(*x*).

In such a case, assuming the provability of $(Ax)F(x)$ must yield a contradiction from which we can either conclude the inconsistency of **PA**, or deduce the unprovability of $(Ax)F(x)$ in a simply consistent **PA**.

The latter conclusion could also then imply the Platonistic existence of some model of **PA** in which $(Ax)\sim PRF'(r, x)$ does not hold true, even though $(An)\sim PRF'(r, n)$ is true in every model of **PA**.

Clearly, in this case, the *omega*-provability of arithmetical assertions is better viewed as a reflection of the 'vagueness or ill-definedness' of some predicate inherent in their construction, where the range of values for which the predicate is satisfied under interpretation does not form a well-defined totality that can be validly referenced within a logical proposition.

### 1.16. The Platonistic axiom of reducibility and PA

A second possibility (*not unrelated to the first*) is that the formal system itself may have been defined Platonistically.

Thus the meta-equivalence in standard first order Arithmetic PA, that:

> **PA** proves: $F(x)$   <=>   **PA** proves: $(Ax)F(x)$

can also be viewed as reflecting elements of the Platonistic Axiom of Reducibility of the Principia (*or the Comprehension Axiom of set theory : Gödel 1931 p178*).

It essentially meta-postulates that the closure of every provable well-formed formula of **PA** containing free variables yields a provable well-formed formula of **PA**.

In every interpretation, the above can be taken to assert that if, for a well-formed formula $F(x)$, $x$ satisfies $F(x)$ whenever $x$ is representable in the formal system, then we can conclude $(Ax)F(x)$ in the interpretation, and hence that the domain of $x$ is well-defined (*i.e. it forms a well-defined set in set theory*) under the interpretation.[5]

**1.17. What is the significance of: $(Ax)F(x)$ is *omega*-provable?**

To sum up, we note that, if $(Ax)F(x)$ is *omega*-provable, then the following hold.

(*i*)   '~(**PA** proves: $(Ax)F(x)$)' **holds meta-mathematically**.

This can be taken to intuitively imply either that there is some model of **PA** in which there is a non-constructive element (*interpreted Platonistically*) in the domain of $x$ that is not representable in **PA**, and is such that it does not satisfy the predicate $F(x)$, or that the totality of values for which the predicate '$F$' holds over the domain of $x$ in some model of **PA** is ill-defined.

(*ii*)   '$(An)$(**PA** proves: $F(n)$)' **holds meta-mathematically, where the domain of $n$ is the set $N$ of natural numbers**.

This intuitively implies that, in every model of **PA**, the predicate $F(x)$ is satisfied constructively by every interpretation of the numeral $n$ in the domain of $x$.

---

[5]  This indicates that we may need to formally distinguish between the meta-statement '$F(x)$ is provable', and the meta-statement '$(x)F(x)$ is provable'. The former may reasonably be taken to intuitively assert that a well-formed formula '$F(x)$' is provable for any given, even if unspecified, element $x$ of every domain in every model. The latter, however, is intuitively the stronger assertion that '$F(x)$' is provable for every element $x$ of any domain; it implicitly postulates that the domain of $x$, in every model, is a well-defined totality for the predicate '$F$', and so it can be validly referenced in a logical proposition. Equating the two can thus be viewed as unreasonable, and so counter-intuitive.

(*iii*)  '(A*x*)(**PA** proves: *F*(*x*))' **holds** meta-**mathematically**.

> This includes (*ii*), but intuitively implies additionally that if, in any model of **PA**, there is (*interpreted Platonistically*) some element in the domain of *x* that is not an interpretation of the numeral ***n***, but is otherwise representable in **PA**, then it necessarily satisfies *F*(*x*).

**1.18.  Why is first order Arithmetic not *omega*-constructive?**

The choice of our formal system of first order Arithmetic may thus depend on whether it is our (*tacit*) intent to deny validity to well-formed formulas that reflect ill-defined elements under every interpretation, or to admit them as well-formed formulas and accept the intuitively artificial conclusion that they are formally unprovable, but translate into true propositions under interpretation.

The question arises: Assuming **PA** corresponds to the latter alternative, can we suitably modify **PA** so that, similar to the development of standard set theory ZF, we arrive at a **PA** system where we prevent the representation of predicates with ill-defined totalities of the values that satisfy them under interpretation, and where we eliminate the non-constructive elements that underlie the logic of **PA**?

In the next section, we therefore address the question: Why is standard first order Arithmetic **PA** not *omega*-constructive?

## Chapter 2.  Why is first order Arithmetic not *omega*-constructive?

**Contents**

2.0.  Introduction
2.1.  The Representation Theorem
2.2.  What is the formal representation of *f*(*x*, *y*) in **PA**?

**2.0. Introduction**

In this second section we locate specific non-constructivity in **PA** and trace it to the particular *closed Induction* axiom schema of **PA**, and the *strong Generalisation* rule of inference of first order predicate calculus.

**2.1. The Representation Theorem**

We start by noting that, in the argument that a simply consistent **PA** cannot be *omega-*constructive, the only overt assumptions (*cf. Gödel 1931 p190-191*) were that the axioms and rules of inference of **PA** are recursively definable, and Podnieks Representation Theorem (Section 3.3).

The constructive nature of the first assertion is easily established from first principles by direct examination of the 45 functions and relations defined recursively by Gödel (*cf. Gödel 1931 p182-186*).

We thus take a closer look at the second assertion that every Turing-computable function $f(x, y)$ can be represented by a **PA** formula $F(x, y, z)$ such that, for every set of natural numbers $k, m, n$, if $f(k, m) = n$ then:

(i)     **PA** proves: $F(\mathbf{k}, \mathbf{m}, \mathbf{n})$.

(ii)    **PA** proves: $(Az)(\sim(z=\mathbf{n}) \Rightarrow \sim F(\mathbf{k}, \mathbf{m}, z))$.

## 2.2. What is the formal representation of $f(x, y)$ in PA?

Assuming that a function is Turing-computable if and only if it is recursive, we consider the case where $f(x, y)$ is a recursive function, defined by:

(i)    $f(x, 0) = g(x)$,

(ii)   $f(x, (y+1)) = h(x, y, f(x, y))$,

(iii) where $g(x_1)$ and $h(x_1, x_2, x_3)$ are recursive functions of lower rank[6] that are represented in **PA** by well-formed formulas $G(x_1, x_2)$ and $H(x_1, x_2, x_3, x_4)$.

Now Mendelson shows (*p131 & p132*) that $f(x_1, x_2)$ is represented by the well-formed formula $F(x_1, x_2, x_3)$:

$(Eu)(Ev)[((Ew)(Bt(u, v, \mathbf{0}, w)\ \&\ G(x_1, w)))\ \&\ Bt(u, v, x_2, x_3)\ \&\ (Aw)(w < x_2 \Rightarrow (Ey)(Ez)(Bt(u, v, w, y)\ \&\ Bt(u, v, (w+\mathbf{1}), z)\ \&\ H(x_1, w, y, z))]$.

---

[6] In the web translation of *Gödel 1931 Proposition V* this is referred to as the 'degree' of the recursive function.

## 2.3. Gödel's Beta function

Defining Gödel's Beta-function[7], $\beta(x_1, x_2, x_3) = rm(1+(x_3+1).x_2, x_1)$, where $rm(x, y)$ is the remainder obtained on dividing $y$ by $x$, the well-formed formula $Bt(x_1, x_2, x_3, x_4)$ in §**2.2** above is a representation in **PA** of $\beta(x_1, x_2, x_3)$, and is the well-formed formula (*Mendelson p131*):

$$(Ew)( x_1 = ((1 + (x_3 + 1). x_2).w + x_4) \ \& \ (x_4 < 1 + (x_3 + 1). x_2)).$$

## 2.4. Can $F(x, y, z)$ create a logical paradox in PA?

Having represented $f(x, y)$ in **PA** by $F(x, y, z)$ for every set of natural numbers $k, m, n$ such that $f(k, m) = n$, the question arises: Is the well-formed formula $F(x, y, z)$ well-defined in **PA**?

In other words, can we ensure that the recursive predicate $f(x, y)$ has not been represented by tacitly postulating a well-formed formula $F(x, y, z)$ of **PA** as well-defined, which may implicitly admit the possibility of a formal logical paradox in **PA**?

## 2.5. The '*Liar*' paradox

To see the significance of this question, we review the earlier paradox, loosely ascribed to Epimenides, which arises obviously if we introduce by definition a '*Liar*' expression, which reads as "The '*Liar*' proposition is a lie", as a valid proposition within the language in which the expression is constructed.

---

[7] Gödel's Beta-function (*Mendelson p131; Gödel 1931 Lemma 1*) may be loosely taken to represent a given finite sequence. Formally, for any given sequence $f(1), f(2), ..., f(n)$, we can always construct some $x_1$ and $x_2$ such that, for all $i<n, \beta(x_1, x_2, i) = f(i)$. We note, however, that the denumerable sequences $f(1), f(2), ..., f(n), m_k$ where $k>0$ and $m_k?m_l$ if $k?l$, are represented by denumerable, distinctly different, Beta-functions $\beta(x_{k,1}, x_{k,2}, i)$ respectively. There are thus denumerable pairs $x_{k,1}, x_{k,2}$ for which $\beta(x_{k,1}, x_{k,2}, i)$ represents any given sequence $f(1), f(2), ..., f(n)$.

Clearly the grammar of an ordinary language, which governs the construction of words and propositions of the language for use as a means of general unrestricted communication, does not contain the specification necessary to prohibit the creation by definition of such vague, or ill-defined, self-referential expressions that have the formal form of a proposition, but which are devoid of both communicative content and meaning.

**2.6. The '*Russell*' paradox**

Russell's expression, within set theory, for a set '*Russell*' whose elements are all, and only, those sets of the Theory that do not belong to themselves, is equally paradoxical.

Here also the axioms and logical rules of inference of the theory, which govern the construction of sets for use in a more restricted and precise language of mathematical communication, do not contain the specification necessary to satisfactorily prohibit the creation by definition of vague, ill-defined or inconsistent entities that apparently harbour concealed concepts involving unruly infinite sets.

**2.7. Creation through definition**

We note that definitions are essentially arbitrary assignments of convenient names within a language or theory for expressing complex reasoning in a compact, more readable, and easier to grasp form.

The paradoxes indicate that unrestricted definitions may admit creation of well-formed expressions within the language or theory that lead to an inconsistency. Clearly such expressions must be considered as ill-defined within the language or theory.

It follows that a language or theory that admits such well-formed expressions as well-defined cannot be reasonably assumed to be simply consistent.

**2.8. Well-defined well-formed formulas of PA**

Now the formal response towards the containment (*as distinct from the resolution*) of such paradoxes in **PA** appears to be the *implicit* thesis that if a well-formed formula is provable in **PA**, then it necessarily represents a function or predicate that is well-defined in **PA**.

We may, therefore, reasonably assume that such well-formed formulas do not affect the formal consistency of the axioms of **PA**.

Thus, to ensure that $F(x_1, x_2, x_3)$ does not introduce an inconsistency in **PA**, Mendelson argues (*p134*) that $F(x_1, x_2, x_3)$ strongly represents $f(x_1, x_2)$ in **PA**, in the sense that:

    **PA** proves: $(E_1x_3)F(x_1, x_2, x_3)$,

where '$E_1x_3$' defines uniqueness (*Mendelson p79*).

However, in the light of our argument in the preceding paragraph, if $F(x_1, x_2, x_3)$ admits construction of well-formed formulas in **PA** that translate, under interpretation, as circular, ill-defined or vaguely defined predicates such as those that lead to the logical antinomies, then the strong representability of $F(x_1, x_2, x_3)$ in **PA** provides reasonable grounds for arguing that **PA** is itself an ill-defined system.

**2.9. What does PA prove constructively by $(E_1x_3)F(k, m, x_3)$?**

We therefore consider in some detail the significance of the formal provability, in **PA**, of the well-formed formula $(E_1x_3)F(k, m, x_3)$ for any given numerals $k, m$:

**PA** proves: $(E_1x_3)(Eu)(Ev)[((Ew)(\boldsymbol{Bt}(u, v, \boldsymbol{0}, w) \And \boldsymbol{G}(k, w))) \And \boldsymbol{Bt}(u, v, m, x_3) \And$
$(Aw)(w< m => (Ey)(Ez)(\boldsymbol{Bt}(u, v, w, y) \And \boldsymbol{Bt}(u, v, (w+\boldsymbol{1}), z) \And \boldsymbol{H}(k, w, y, z))].$

Now "**PA** proves: $(E_1x_3)\boldsymbol{F}(k, m, x_3)$" intuitively asserts that, for any given natural numbers $k, m$, we can *construct* natural numbers $x_3$ (*uniquely*), $u, v$ that define a Beta-function such that $\beta(u, v, 0) = g'(k)$ and, for all $i<m$, $\beta(u, v, i) = h'(k, i, f'(k, i))$, and $\beta(u, v, m) = x_3$, where $f'(x_1, x_2), g'(x_1)$ and $h'(x_1, x_2, x_3)$ are any recursive functions that are formally represented by $\boldsymbol{F}(x_1, x_2, x_3), \boldsymbol{G}(x_1, x_2)$ and $\boldsymbol{H}(x_1, x_2, x_3, x_4)$ respectively such that:

(i) $f'(k, 0) = g'(k)$,

(ii) $f'(k, (y+1)) = h'(k, y, f'(k, y))$ for all $y<m$,

(iii) $g'(x_1)$ and $h'(x_1, x_2, x_3)$ are recursive functions that are assumed to be of lower rank than $f'(x_1, x_2)$.

For any arbitrarily given sequence of natural numbers $k_0, k_1, \ldots, k_n$, we can clearly determine an infinity of values of $u_p, v_p, k_p$ such that the Beta-function $\beta(u_p, v_p, i) = k_i$ for all $0 =< i =< n$, and $\beta(u_p, v_p, (n+1)) = k_p$.

Hence the **PA** provability of $(E_1x_3)\boldsymbol{F}(k, m, x_3)$ is simply the intuitive assertion that, for any given natural numbers $k$ and $m$, we can always *construct* some (*non-unique*) pair of natural numbers $u, v$ such that the Beta-function $\beta(u, v, i)$ represents the first $m$ terms, i.e. $f'(k, 0), f'(k, 1), \ldots, f'(k, m)$, which are common to every primitive recursive function such as $f'(x_1, x_2)$ that is formally represented by $\boldsymbol{F}(x_1, x_2, x_3)$.

We can see that this is constructively provable for any given natural numbers $k$ and $m$ since, given $F(x_1, x_2, x_3)$, we can determine a suitable Turing machine to *construct* the sequence $f'(k, 0), f'(k, 1), \ldots, f'(k, m)$ uniquely and verify the assertion.

We note, however, that if $g$ is the Gödel-number of any Turing-computable proof of $(E_1x_3)F(k, m, x_3)$ as above then, trivially, $g>m$.

We cannot thus conclude from the above argument that there is necessarily a Turing-computable constructive proof of the existence of a **PA**-formula that yields the non-terminating sequence $f(k, 0), f(k, 1), \ldots, f(k, m), \ldots$, which we could therefore take as representing the function $f(k, y)$ uniquely.

### 2.10. Mendelson's proof of PA proves: $(E_1x_3)F(k, m, x_3)$

The Turing-decidability of the assertion $(E_1x_3)F(k, m, x_3)$ is reflected in Mendelson's proof of the above (*p133*), which relies on an essentially constructive meta-induction argument as below that does not appeal to *strong Generalisation*:

> **PA** proves: $(E_1x_3)F(0, 0, x_3)$
>
> $(Am)[[$**PA** proves: $(E_1x_3)F(0, m, x_3)]$ => [**PA** proves: $(E_1x_3)F(0, m+1, x_3)]]$
>
> $(Am)[$**PA** proves: $(E_1x_3)F(0, m, x_3)]$
>
> $(Ak)(Am)[[$**PA** proves: $(E_1x_3)F(k, m, x_3)]$ => [**PA** proves: $(E_1x_3)F(k+1, m, x_3)]]$
>
> $(Ak)(Am)[$**PA** proves: $(E_1x_3)F(k, m, x_3)]$

## 2.11. What does $(E_1x_3)F(x_1, x_2, x_3)$ assert?

Turning next to the issue of strong representability, we consider the meta-assertion:

**PA** proves: $(E_1x_3)(Eu)(Ev)[((Ew)(\boldsymbol{Bt}(u, v, \boldsymbol{0}, w)\ \&\ \boldsymbol{G}(x_1, w)))\ \&\ \boldsymbol{Bt}(u, v, x_2, x_3)$
$\&\ (Aw)(w < x_2 => (Ey)(Ez)(\boldsymbol{Bt}(u, v, w, y)\ \&\ \boldsymbol{Bt}(u, v, (w+\boldsymbol{1}), z)\ \&$
$\boldsymbol{H}(x_1, w, y, z))].$

Now "**PA** proves: $(E_1x_3)F(x_1, x_2, x_3)$" intuitively asserts that, for variable $x_1, x_2$, we can similarly *construct* natural numbers $x_3$ (*uniquely*), $u, v$ that define a Beta-function such that $\beta(u, v, 0) = g'(x_1)$ and, for all $i =< x_2, \beta(u, v, i) = h'(x_1, i, f'(x_1, i))$, where $f'(x_1, x_2)$, $g'(x_1)$ and $h'(x_1, x_2, x_3)$ are any recursive functions that are formally represented by $F(x_1, x_2, x_3)$, $G(x_1, x_2)$ and $H(x_1, x_2, x_3, x_4)$ respectively such that:

(i) $f'(x_1, 0) = g'(x_1)$,

(ii) $f'(x_1, (y+1)) = h'(x_1, y, f'(x_1, y))$ for all $y < x_2$,

(iii) $g'(x_1)$ and $h'(x_1, x_2, x_3)$ are recursive functions that are assumed to be of lower rank than $f'(x_1, x_2)$.

Hence the **PA** provability of $(E_1x_3)F(x_1, x_2, x_3)$ interprets as the intuitive assertion that, even for variable $x_1$ and $x_2$, we can *determine* some pair of natural numbers $u, v$ such that the Beta-function $\beta(u, v, i)$ represents any indeterminate number, $x_2$, of terms, i.e. the sequence $f'(x_1, 0), f'(x_1, 1), \ldots, f'(x_1, x_2)$, of every primitive recursive function such as $f'(x_1, x_2)$ that is formally represented by $F(x_1, x_2, x_3)$.

## 2.12. The introduction of Platonism into PA

In other words, the **PA** provability of $(E_1x_3)F(x_1, x_2, x_3)$ corresponds to the intuitive, non-constructive, assertion that there is (*interpreted Platonistically*) some Beta-function $\beta(u,$

$v$, $i$) that represents the sequence $f(x_1, 0), f(x_1, 1), \ldots, f(x_1, x_2)$ for any unspecified, and indeterminate, length $x_2$.

The above is, essentially, the assertion that $F(x_1, x_2, x_3)$ uniquely represents $f(x_1, x_2)$ in **PA**. It implicitly implies that the non-terminating sequence of values $f(x_1, 0), f(x_1, 1), \ldots, f(x_1, x_2), \ldots$ forms a well-defined totality, which can be uniquely named and referenced in **PA** by some (*unique up to isomorphism*) Beta-function $\beta_f(u, v, i)$.

The **PA** provability of $(E_1x_3)F(x_1, x_2, x_3)$ thus corresponds to a non-constructive assertion of the Platonistic existence of such a totality.

**2.13. How does PA prove the Platonistic assertion** $(E_1x_3)F(x_1, x_2, x_3)$**?**

The question arises: How, then, does **PA** prove an essentially Platonistic assertion such as $(E_1x_3)F(x_1, x_2, x_3)$?

After establishing the argument by which we can show that **PA** proves: $(E_1x_3)F(\boldsymbol{k}, \boldsymbol{m}, x_3)$, Mendelson notes, indicating outlines of the proof (*p134*), that, for primitive recursive functions (*Mendelson p120*):

> "The part of the proof given above to show that the recursion rule does not lead out of the class of representable functions can be strengthened to yield strong representability ..."

Mendelson's indicated argument appears to be that we can establish by meta-induction:

> **PA** proves: $(E_1x_3)F(\boldsymbol{0}, \boldsymbol{0}, x_3)$
>
> [**PA** proves: $(E_1x_3)F(\boldsymbol{0}, x_2, x_3)$]  =>  [**PA** proves: $(E_1x_3)F(\boldsymbol{0}, x_2+1, x_3)$]
>
> **PA** proves: $(E_1x_3)F(\boldsymbol{0}, x_2, x_3)$ => $(E_1x_3)F(\boldsymbol{0}, x_2+1, x_3)$

**PA** proves: $(Ax_2)((E_1x_3)F(0, x_2, x_3) \Rightarrow (E_1x_3)F(0, x_2+1, x_3))$ ... *Gen*

**PA** proves: $(Ax_2)(E_1x_3)F(0, x_2, x_3)$ ... *Induction*

[**PA** proves: $(Ax_2)(E_1x_3)F(x_1, x_2, x_3)$] $\Rightarrow$ [**PA** proves: $(Ax_2)(E_1x_3)F(x_1+1, x_2, x_3)$]

**PA** proves: $(Ax_2)(E_1x_3)F(x_1, x_2, x_3) \Rightarrow (Ax_2)(E_1x_3)F(x_1+1, x_2, x_3)$

**PA** proves: $(Ax_1)((Ax_2)((E_1x_3)F(x_1, x_2, x_3) \Rightarrow (Ax_2)(E_1x_3)F(x_1+1, x_2, x_3))$

**PA** proves: $(Ax_1)(Ax_2)(E_1x_3)F(x_1, x_2, x_3)$

**PA** proves: $(E_1x_3)F(x_1, x_2, x_3)$

Assuming that this correctly represents Mendelson's intention, the validity of the above reasoning gives the necessary validity to $F(x_1, x_2, x_3)$ as a well-defined well-formed formula that is consistent with the axioms of **PA**.

### 2.14. At which point does Platonism enter PA formally?

The question arises: At which point does **PA** formalise Platonistic elements?

We note that, in the above argument, the non-constructive conclusions in the reasoning are those thart arise from application of the *closed Induction* axiom schema (*Mendelson p103*) and the *strong Generalisation* rule of inference of first order predicate calculus (*Mendelson p57*):

(i)  *Closed Induction* : $(F(0) \& ((Ax)(F(x) \Rightarrow F(x+1)))) \Rightarrow (Ax)(F(x))$

(ii) *Strong Generalisation* : $(Ax_i)F(x_i)$ follows from $F(x_i)$

These are the only assumptions that can be objected to intuitionistically in the above argument to provide legitimacy for $F(x_1, x_2, x_3)$ as a well-defined well-formed formula that is consistent with the axioms of **PA**.

**2.15. Is PA necessarily Platonistic?**

The question arises: Assuming that **PA** is simply consistent, does the above argument, that **PA** is a non-constructive system, imply that **PA** is necessarily Platonistic?

To put the query in perspective, we note that the meta-assertion:

> **PA** proves: $F(x_i)$

intuitively represents a denumerable infinity of formal statements:

> **PA** proves: $F(k)$,

where only the denumerable terms that can be expressed in **PA**, such as $k$, are substituted for $x_i$.

However, the meta-assertion:

> **PA** proves: $(Ax_i)F(x_i)$

is a single formal assertion that, under interpretation in a model of **PA**, may translate into an infinitely-compound proposition that is not necessarily denumerably-compound.

Since the two are formally equivalent in the presence of *strong Generalisation*, the question is then: Does *strong Generalisation* formalise Platonism[8] in **PA**?

---

[8] Gödel's Platonistic interpretation of *strong Generalisation* is implicit when he notes (Gödel *1931 p191 footnote 48a*): "The true reason for the incompleteness which attaches to all formal systems of mathematics lies ... in the fact that the formation of higher and higher types can be continued into the transfinite ..., while, in every formal system, only

**2.16. Is every model of PA necessarily infinite, but not denumerable?**

The above equivalence in **PA** also raises the question: Is every model of **PA** necessarily, infinite but not denumerable, since **PA** is shown to contain formally non-constructive elements.

As we have argued in §**1.17** of the previous section "*Chapter 1. First order Arithmetic is not omega-constructive*", there is an in-built ambiguity in the construction of **PA** once it admits *omega*-provable propositions.

We then face the possibility of predicates giving rise to well-defined well-formed formulas in **PA** that may not faithfully reflect our intuitive understanding of the natural numbers, but that are not obviously transfinite in nature.

This ambiguity, prima facie, does not admit an unequivocal answer to the above question.

**2.17. Is there a constructive PA?**

The question thus arises: Can we formally represent our intuitive Arithmetic of the natural numbers less ambiguously, and more faithfully?

We address this question in the next section: Is there an *omega*-constructive first order Arithmetic?

---

countable many are available. Namely, one can show that the undecidable propositions which have been constructed here always become decidable through adjunction of sufficiently high types (*e.g. of the type **omega** to the system **P***). A similar result holds for axioms of set theory."

# Chapter 3.  Is there an *omega*-constructive first order Arithmetic?

**Contents**

## 3.0.  Introduction

In this third section we formally construct and briefly examine various systems of first order Arithmetic. We argue that we cannot use Gödel's reasoning to prove that every primitive recursive function (*Mendelson p120*) can be strongly represented formally in all such systems. Hence it does not universally establish the existence of undecidable propositions that are true under interpretation.

We finally construct an *omega*-constructive system of first order Arithmetic, *omega*-**PA**, where we replace the non-constructive *closed Induction* axiom schema by a constructive, *open Induction* axiom schema.

We base *omega*-**PA** on a modified first order predicate calculus where we replace the *strong Generalisation* rule of inference I(*ii*) by a *weaker omega-Constructivity* rule of inference I(*iii*).

We show that every model of *standard* **PA** is a model of *omega*-**PA**. We show further that every primitive recursive function cannot be strongly represented formally in *omega*-**PA** if it is simply consistent. So, in *omega*-**PA**, Gödel's reasoning does not establish the existence of undecidable propositions that are true under interpretation.

We argue in this, and the next, section that Gödel mistakenly concluded (*Gödel p190-191*) that the existence of undecidable propositions that are true under interpretation must follow from his reasoning in all the above formal systems of Arithmetic.

We begin by reviewing in detail Mendelson's definition (*p103*) of the system **PA**. This is essentially the formal system in which Gödel constructed his undecidable propositions.

### 3.1. The primitive symbols of a general first order predicate calculus K

The system **PA** is based on a particularisation of the general first order predicate calculus **K** whose primitive symbols are (*Mendelson p57*):

'~' (*not*), '=>' (*implies*), '(' (*left parenthesis*), ')' (*right parenthesis*), ',' (*comma*), '$x_1, x_2, ...$' (*denumerably many individual variables*), '$A_{(n, j)}, (n, j >= 1)$' (*denumerably many predicate letters*), '$f_{(n, j)}, (n, j >= 1)$' (*denumerably many function letters*), and '$a_i, (i >= 1)$' (*denumerably many constants*).

### 3.2. The logical axioms of K

The logical axioms of **K** are:

    **K**(1)   $A \Rightarrow (B \Rightarrow A)$

**K**(2)  $(A \Rightarrow (B \Rightarrow C)) \Rightarrow ((A \Rightarrow B) \Rightarrow (A \Rightarrow C))$

**K**(3)  $(\sim B \Rightarrow \sim A) \Rightarrow ((\sim B \Rightarrow A) \Rightarrow B)$

**K**(4)  $(Ax_i)A(x_i) \Rightarrow A(t)$, if $A(x_i)$ is a well-formed formula of **K** and $t$ is a term of **K** free for $x_i$ in $A(x_i)$.

**K**(5)  $(Ax_i)(A \Rightarrow B) \Rightarrow (A \Rightarrow (A \Rightarrow (Ax_i)B))$ if $A$ is a well-formed formula of **K** containing no free occurences of $x_i$.

### 3.3. The rules of inference of K

The rules of inference of **K** are:

**I**(*i*)  *Modus Ponens*: **B** follows from $A$ and $A \Rightarrow B$.

**I**(*ii*)  *Strong Generalisation*: $(Ax_i)A$ follows from $A$.

### 3.4. Gödel's formal system PA

The particular system of first order predicate calculus that Mendelson uses (*p103*), to develop Gödel's reasoning in a formal first order Arithmetic, contains:

**A**(*i*)  a single predicate letter '**A**' to represent equality $A(t, s)$ between the terms $t$ and $s$, which we express as '$t=s$';

**A**(*ii*)  one individual constant '$a$' which we write as '0';

**A**(*iii*)  three function letters $f$, $g$ and $h$ to represent, respectively, the successor of '$t$', which we write as '$t'$', addition, which we write as '$t+s$', and multiplication, which we write as '$t.s$'.

The axioms of Arithmetic that Gödel considered were Dedekind's axioms for number theory, known as Peano's Postulates, which are essentially as below (*not necessarily independent*):

**A(1)**   $(x_1 = x_2) \Rightarrow ((x_1 = x_3) \Rightarrow (x_2 = x_3))$

**A(2)**   $(x_1 = x_2) \Rightarrow (x_1' = x_2')$

**A(3)**   $0 ? x_1'$

**A(4)**   $(x_1' = x_2') \Rightarrow (x_1 = x_2)$

**A(5)**   $(x_1 + 0) = 0$

**A(6)**   $(x_1 + x_2') = (x_1 + x_2)'$

**A(7)**   $x_1.0 = 0$

**A(8)**   $x_1.x_2' = (x_1.x_2) + x_1$

**A(9)**   $(x_1 ? 0) \Rightarrow (Ex_2)(x_1 = x_2')$

and the two meta-axiom schemas (*more accurately, rules of inference*) of *closed* and *open Induction*:

**A(10)**   For any wff $F(x)$ of **PA**, from $(F(0)$ and $(Ax)(F(x)) \Rightarrow F(x'))$, we may conclude $(Ax)F(x)$

**A(11)**   For any wff $F(x)$ of **PA**, from $(F(0)$ and $(F(x)) \Rightarrow F(x'))$ we may conclude $F(x)$

### 3.5. General Arithmetics : *weak*-GA

If we eliminate *strong Generalisation* **I**(*ii*), *closed Induction* **A**(10) and *open Induction* **A**(11) from **PA**, we obtain a weak first order general Arithmetic, *weak*-**GA**.

*Weak*-**GA** is thus the axioms **A**(1) to **A**(9) of first order Arithmetic, as specified in §**3.4**, added to a first order predicate calculus **A**(*i*) to **A**(*iii*), as specified in §**3.4**, with only the *Modus Ponens* rule of inference **I**(*i*).

Since all axioms of *weak*-**GA** are provable in **PA**, every interpretation and model of **PA** is an interpretation and model of *weak*-**GA**.

However, in the absence of *strong Generalisation* **I**(*ii*), and *closed Induction* **A**(10), the argument of §**2.13**, establishing the *strong representability* of every recursive function, cannot be carried out in *weak*-**GA**.

So, in *weak*-**GA**, Gödel's reasoning does not establish the existence of undecidable propositions that are true under interpretation.

### 3.6. General Arithmetics : *strong*-GA

If we strengthen *weak*-**GA** by adding the *strong Generalisation* rule of inference **I**(*ii*), we arrive at a stronger first order general Arithmetic, *strong*-**GA**.

*Strong*-**GA** is thus the axioms **A**(1) to **A**(9) of first order Arithmetic, as specified in §**3.4**, added to a first order predicate calculus **A**(*i*) to **A**(*iii*), as specified in §**3.4**, with both the *Modus Ponens* rule of inference **I**(*i*) and the *strong Generalisation* rule of inference **I**(*ii*).

Since all axioms of *strong*-**GA** are provable in **PA**, every interpretation and model of **PA** is an interpretation and model of *strong*-**GA**.

However, in the absence of *closed Induction* **A**(10), the argument of §**2.13**, establishing the strong representability of every recursive function, still cannot be carried out in *strong*-**GA**.

So, in *strong*-**GA**, Gödel's reasoning does not establish the existence of undecidable propositions that are true under interpretation.

**3.7. Peano Arithmetics :** *weak*-**PA**

If we add *closed Induction* **A**(10) to *weak*-**GA**, we arrive at a weak first order Peano Arithmetic, *weak*-**PA**.

*Weak*-**PA** is thus the axioms **A**(1) to **A**(10) of first order Arithmetic, as specified in §**3.4**, added to a first order predicate calculus **A**(*i*) to **A**(*iii*), as specified in §**3.4**, with only the *Modus Ponens* rule of inference **I**(*i*).

Since all axioms of *weak*-**PA** are provable in **PA**, every interpretation and model of **PA** is an interpretation and model of *weak*-**PA**.

However, in the absence of *strong Generalisation* **I**(*ii*), the argument of §**2.13**, which establishes the strong representability of every recursive function, cannot be carried out in *weak*-**PA**.

Hence, in *weak*-**PA** also, Gödel's reasoning does not establish the existence of undecidable propositions that are true under interpretation.

We note also that it is not obvious whether *weak*-**PA** proves *weak-Induction* **A**(11).

**3.8. Peano Arithmetics :** *Standard* **PA**

We arrive at the standard first order Peano Arithmetic, *standard* **PA**, used by Gödel, when we add *strong Generalisation* **I**(*ii*) to *weak*-**PA**.

*Standard* **PA** is thus the axioms **A**(1) to **A**(10) of first order Arithmetic, as specified in §**3.4**, added to a first order predicate calculus **A**(*i*) to **A**(*iii*), as specified in §**3.4**, with both the rules of inference, *Modus Ponens* **I**(*i*), and *strong Generalisation* **I**(*ii*).

Since *standard* **PA** establishes the *strong* representability of every recursive function, we are now able to formally derive Gödel's undecidable propositions in *standard* **PA**.

However, as we have argued in the earlier sections, *strong* representability now admits closed well-formed formulas that refer to totalities, under interpretation, which are not reflected in the intuitive Arithmetic that is intended to be formalised by *standard* **PA**.

**3.9.** *Closed Induction* **implies** *open Induction* **in PA**

We note that, in *standard* **PA**, *closed Induction* **A**(10) implies *open Induction* **A**(11).

Thus, from the *closed Induction* axiom schema **A**(10) in *standard* **PA**:

>   *Standard* **PA** proves: (*F(0)* **&** *Standard* **PA** proves: (A*x*)(*F(x)* => *F(x+1)*))

>   => *Standard* **PA** proves: (A*x*)*F(x)*

we have that if:

>   *Standard* **PA** proves: *F(0)* **&** *Standard* **PA** proves: *F(x)* => *F(x+1)*

then:

>   *Standard* **PA** proves: *F(x)* => *F(x+1)*

>   *Standard* **PA** proves: (A*x*)(*F(x)* => *F(x+1)*)    ... Gen

>   *Standard* **PA** proves: *F(0)* **&** *Standard* **PA** proves: (A*x*)(*F(x)* => *F(x+1)*)

>   *Standard* **PA** proves: (A*x*)*F(x)*

*Standard* **PA** proves: *F*(*x*)

We thus have that:

*Standard* **PA** proves: *F*(**0**)  &  *Standard* **PA** proves: (*F*(*x*) => *F*(*x*+**1**)))

=>  *Standard* **PA** proves: *F*(*x*)

Hence *standard* **PA** proves the *open Induction* axiom schema **A**(11), and so **PA** and *standard* **PA** are identical.

### 3.10. Can PA admit models of transfinite ordinals?

Now we note that:

**PA** proves: *x*+1 = 1+*x*

**PA** proves: (A*x*)(*x*+1 = 1+*x*).

It follows that Cantor's transfinite *ordinal*-Arithmetic **CA**, in which *omega* is a term for which we have that "*omega* +1" is not equal to "1+ *omega*", is not a model of **PA**.

### 3.11. *omega*-constructive first order Arithmetics

We next define a range of *omega*-constructive first order systems of Arithmetic, by introducing an intuitionistically unobjectionable *omega-Constructivity* rule of inference:

**I**(*iii*)   *omega-Constructivity*: *F*(*x*) is equivalent to (A*x*)["*x* is a numeral" => *F*(*x*)].

Since "*x* is a numeral" is recursively definable, we note that **I**(*iii*) is a recursively definable rule of inference that is equivalent to the concept of *omega-constructivity* introduced in "*Chapter 1. First order Arithmetic is not omega-constructive*".

### 3.12. *omega*-constructive Arithmetic : *omega*-GA

If we strengthen *weak*-**GA** by adding the *omega-Constructivity* rule of inference **I**(*iii*), we arrive at a formal system of constructive *omega*-**GA**.

*omega*-**GA** is thus the axioms **A**(1) to **A**(9) of first order Arithmetic, as specified in §**3.4**, added to a first order predicate calculus **A**(*i*) to **A**(*iii*), as specified in §**3.4**, along with the *Modus Ponens* rule of inference **I**(*i*) and the *omega-Constructivity* rule of inference **I**(*iii*).

Since all axioms of *omega*-**GA** are provable in **PA**, every interpretation and model of **PA** is an interpretation and model of *omega*-**GA**.

Also, in the absence of *strong Generalisation* **I**(*ii*), and *closed Induction* **A**(10), the argument of §**2.13**, establishing the *strong* representability of every recursive function, cannot be carried out in *omega*-**GA**.

However, as we have shown in the earlier section, "*Chapter 1. First order Arithmetic is not omega-constructive*", in any system of Arithmetic with the *omega-Constructivity* rule of inference **I**(*iii*), the assumption that every recursive function can be *strongly* represented leads to inconsistency in the system.

Hence, in *omega*-**GA**, Gödel's reasoning does not establish the existence of undecidable propositions that are true under interpretation.

### 3.13. The *omega*-constructive Peano Arithmetic : *omega*-PA

If we now add *open Induction* **A**(11) to *omega*-**GA**, we arrive at a constructive formal system of Peano's Arithmetic, *omega*-**PA**.

*omega*-**PA** is thus axioms **A**(1) to **A**(9), **A**(11) of first order Arithmetic, as specified in §**3.4**, added to a first order predicate calculus **A**(*i*) to **A**(*iii*), as specified in §**3.4**, along with the *Modus Ponens* rule of inference **I**(*i*) and the *omega-Constructivity* rule of inference **I**(*iii*).

Since the assumption that every recursive function can be *strongly* represented in any system of Arithmetic with the *omega-Constructivity* rule of inference **I**(*iii*) leads to inconsistency in the system, it follows that Gödel's reasoning does not hold in any interpretation of *omega*-**PA**.

We note that the *open Induction* axiom schema **A**(11) permits us to constructively establish in *omega*-**PA**, for open well-formed formulas *F*(*x*) for which we can prove both the well-formed formulas *F(0)* and *F(x) => F(x+1)*, that the well-formed formula *F*(*x*) is provable, and so necessarily well-defined.

### 3.14. Every model of PA is a model of *omega*-PA

Since all axioms of *omega*-**PA** are provable in **PA**, every interpretation and model of **PA** is an interpretation and model of *omega*-**PA**.

Clearly, even though *omega*-**PA** does not prove *closed Induction* **A**(10), it is, intuitively, a natural and constructive formalisation of Peano's Arithmetic.

Moreover, Gödel's reasoning does not establish the existence of undecidable propositions in *omega*-**PA** in a constructive and intuitionistically unobjectionable way, as he intended (*Gödel 1931 p189*).

### 3.15. The essential difference between PA and *omega*-PA

We note that the essential difference between **PA** and *omega*-**PA** lies simply in their respective meta-assertions about the nature and scope of the valid consequences[9] amongst the elements of the domain of each system under interpretation.

Thus the *strong Generalisation* rule of inference, for instance, asserts that, if $F(x)$ is provable, then $(x)F(x)$ is provable, and so, under every interpretation, the range of values for which the predicate '$F$' is satisfied can be referred to as a well-defined totality.

In other words, the predicate '$F$' is asserted, Platonistically, to be a characteristic of every element of the domain under every interpretation, irrespective of whether the element corresponds to any constructive representation within the formal system **PA** or not.

In contrast, the meta-rule of *omega-Constructivity* rule of inference is the *weaker*, constructive, and intuitionistically unobjectionable assertion that, if $F(x)$ is provable, this only implies that the predicate '$F$' is asserted to be a provable characteristic of every element of the domain under interpretation that has a constructive representations within the formal system *omega*-**PA**.

### 3.16. Can *omega*-PA admit models of transfinite ordinals?

We note that Cantor's transfinite *ordinal*-Arithmetic, **CA**, contains the system of natural numbers and so, in the absence of *strong Generalisation*, all the axioms of *omega*-**PA** translate into true statements, in **CA**, about only the sub-domain **N** of natural numbers.

Thus **CA** is a model of *omega*-**PA**.

---

[9] See also Peregrin - Section 11, Conclusions.

This reflects the fact that *omega*-**PA** permits transfinite interpretations in which relations that hold in any denumerable sub- domain of the interpretation may not hold in the entire domain of the interpretation.

For instance, although:

> *omega*-**PA** proves: $x+1 = 1+x$

we also have:

> ~(*omega*-**PA** proves: $(Ax)(x+1 = 1+x)$).

**3.17. Can we strengthen *omega*-PA**

This raises the natural question: At which point does *omega*-**PA** loose both its non-constructive character, and its power to permit transfinite interpretations as above?

In other words, can we strengthen *omega*-**PA**, and extend the range of its provable well-formed formulas under interpretation, without losing its essentially constructive, and intuitionistically unobjectionable character?

In the next section, we address this question.

# Chapter 4. Are there stronger *omega*-constructive systems of first order Arithmetic?

**Contents**

**4.0.  Introduction**

In this fourth section we begin by briefly reviewing the main arguments of the earlier sections.

We then consider whether, and how, we may strengthen *omega*-**PA**.

We use Parikh's form of the Kreisel-Parson's conjecture to constructively qualify quantification. This now allows a constructive formal system to refer to selective properties of non-constructive elements under interpretation

We finally construct a *strong omega$_2$*-**PA**, and show that the axioms, interpretations and models of *strong omega$_2$*-**PA**, and the axioms, interpretations and models of **PA**, are identical, but have significantly different consequences.

We conclude that the formal undecidability of Arithmetical propositions that are true under interpretation is a characteristic not of any relations that are Platonistically inherent in any Arithmetic of the natural numbers, but of the particular formalisation chosen to represent the Arithmetic.

**4.1.  An *omega*-constructive first order Arithmetic**

We argued in §**1** and §**2** that an *omega*-constructive system of Arithmetic is inconsistent with standard first order predicate calculus **PA** with *strong Generalisation*.

### 4.2. *omega*-PA

In §**3** we defined *omega*-**PA**, where we eliminated the *closed Induction* axiom schema **A**(10) of **PA**, and replaced the *strong Generalisation* rule of inference **I**(*ii*) of **PA** by the *weaker omega-Constructivity* rule of inference **I**(*iii*).

### 4.3. Some features of *omega*-PA

We then argued that every model of **PA** is a model of *omega*-**PA**.

We argued further that every primitive recursive function cannot be *strongly* represented in *omega*-**PA**.

We also argued that *omega*-**PA** is intuitively a more natural and constructive formal system of Peano's Arithmetic than **PA**.

We argued that Gödel's reasoning does not establish the existence of undecidable propositions in *omega*-**PA**, which are true under interpretation, in a constructive and intuitionistically unobjectionable way, as he intended (*Gödel p189*).

### 4.4. The difference between PA and *omega*-PA

We argued in §**3.15** that the essential difference between **PA** and *omega*-**PA** lay in their respective meta-assertions about the nature and scope of the valid consequences amongst the elements of the domain of each system under interpretation.

We then posed the question: Can we strengthen *omega*-**PA** so that every model of a stronger *omega*-**PA** is a model of **PA**?

### 4.5. *omega₁*-**PA**

If we now add the *closed Induction* axiom schema **A**(10) to *omega*-**PA**, we arrive at a stronger, yet still constructive, formal system of Peano Arithmetic, *omega₁*-**PA**.

*omega₁*-**PA** is thus axioms **A**(1) to **A**(11) of first order Arithmetic, as specified in §**3.4**, added to a first order predicate calculus **A**(*i*) to **A**(*iii*), as specified in §**3.4**, along with the *Modus Ponens* rule of inference **I**(*i*) and the *omega-Constructivity* rule of inference **I**(*iii*).

Thus **PA** and *omega₁*-**PA** have identical axioms, with identical interpretations and models. Prima facie, though, *omega₁*-**PA** is an intuitively more natural and constructive formalisation of Peano's Arithmetic than **PA**.

Since the assumption that every recursive function can be *strongly* represented in any system of Arithmetic with the *omega-Constructivity* rule of inference **I**(*iii*) leads to inconsistency in the system, it follows that Gödel's reasoning does not hold in any interpretation of *omega₁*-**PA**.

Further, although *omega₁*-**PA** now proves the *closed Induction* axiom schema **A**(11), in the absence of the *strong Generalisation* rule of inference **I**(*ii*) of **PA**, we still cannot conclude either of:

  *omega₁*-**PA** proves: '*Closed Induction*'  =>  '*Open Induction*' holds in *omega₁*-**PA**

  '*Open Induction*' holds in *omega₁*-**PA**  =>  *omega₁*-**PA** proves: '*Closed Induction*'

whereas in **PA** we have both:

  **PA** proves: '*Closed Induction*'  =>  '*Open Induction*' holds in **PA**

'*Open Induction*' holds in **PA**   =>   **PA** proves: '*Closed Induction*'

### 4.6. The conflict with model theory

Now, although **PA** and *omega*$_1$-**PA** have identical axioms, with identical interpretations and models, the difference in their rules of inference yield differing consequences amongst the well-formed formulas of the two formal systems, and hence amongst corresponding statements under interpretation.

The difference is reflected in the fact that Cantor's transfinite ordinal Arithmetic **CA** is not a model of **PA**, since:

> **PA** proves: $x+1 = 1+x$
>
> **PA** proves: $(Ax)(x+1 = 1+x)$.

However, in *omega*$_1$-**PA**, we have that, as in *omega*-**PA**:

> *omega*$_1$-**PA** proves: $x+1 = 1+x$
>
> ~(*omega*$_1$-**PA** proves: $(Ax)(x+1 = 1+x)$).

Further, in the absence of *strong Generalisation*, all the axioms of *omega*$_1$-**PA** translate into true statements, in Cantor's transfinite *ordinal*-Arithmetic **CA**, about only the sub-domain **N** of natural numbers.

Thus **CA** is a model of *omega*$_1$-**PA**.

This appears to conflict with the implicit Platonistic thesis that it is the axioms of a formal system that essentially determine the logical consequences in the interpretation.

### 4.7. Introducing the *omega-Specification* rule of inference

Can we strengthen *omega₁*-**PA** further?

We consider introducing Parikh's form of the Kreisel-Parson's conjecture in *omega₁*-**PA** as a means for specifying essentially the conditions under which the closure of a provable, open, well-formed formula in *omega₁*-**PA** yields a provable, closed, well-formed formula of *omega₁*-**PA**:

> **I**(*iv*)   *omega-Specification*: $(Ax)F(x)$ is provable in **PA** if, for each natural number $n$, there is a proof in **PA** of $F(n)$, and its Gödel-number is less than some natural number $k$ that is independent of $n$.

### 4.8. A yet stronger *omega*-calculus : *omega₂*-**PA**

We thus arrive at a yet stronger, constructive, formulation of the first order predicate calculus, *omega₂*-calculus.

*omega₂*-**PA** is now the axioms **A**(1) to **A**(11) of first order Arithmetic, as specified in §**3.4**, added to a first order predicate calculus **A**(*i*) to **A**(*iii*), as specified in §**3.4**, along with the *Modus Ponens* rule of inference **I**(*i*), the *omega-Constructivity* rule of inference **I**(*iii*) and the *omega-Specification* rule of inference **I**(*iv*).

Thus **PA** and *omega₁*-**PA**, *omega₂*-**PA** all have identical axioms, with identical interpretations and models. Again, *omega₂*-**PA** appears intuitively a more natural and constructive formalisation of Peano's Arithmetic than **PA**.

We note that, as in **PA**, we now have that:

> *omega₂*-**PA** proves: $x+1 = 1+x$

$omega_2$-**PA** proves: $(Ax)(x+1 = 1+x)$

So, like **PA**, $omega_2$-**PA** too does not admit Cantor's transfinite ordinals in any model.

Thus, though $omega_2$-**PA** and $omega_1$-**PA** are both constructive and have identical axioms, these too have significantly different logical consequences, highlighting again the apparent conflict with the implicit Platonistic thesis that it is the axioms of a formal system that essentially determine the logical consequences.

### 4.9. The *Generalisation* rule of inference holds for $omega_2$-**PA**

We note that, in $omega_2$-**PA**, as in **PA**, the *strong Generalisation* rule of inference **I**(*ii*) holds, since we now have the meta-equivalence:

$omega_2$-**PA** proves: $F(x)$ <=> $omega_2$-**PA** proves: $(Ax)F(x)$

It follows that, in $omega_2$-**PA**, as in **PA**, we have:

$omega_2$-**PA** proves: '*Closed Induction*' <=> '*Open Induction*' holds in $omega_2$-**PA**

### 4.10. We cannot infer Gödel's conclusions from his reasoning in $omega_2$-**PA**

However, $omega_2$-**PA** retains the constructive and intuitionistically unobjectionable nature inherited from $omega$-**PA**.

Thus, in sharp contrast to **PA**, $omega_2$-**PA** does not admit the *strong* representation of every recursive function, and so Gödel's reasoning establishing the existence of undecidable propositions cannot be carried out in $omega_2$-**PA**.

This follows from the fact that, if $f(x, y)$ is a recursive function that is represented by a $omega_2$-**PA** well-formed formula $F(x, y, z)$, then the formal provability of $(E_1x_3)F(k, m,$

$x_3$) is the intuitive assertion (*cf.* §**2.9**) that, for any given natural numbers $k$ and $m$, we can always *construct* some (*non-unique*) pair of natural numbers $u$, $v$ such that the Beta-function $\beta(u, v, i)$ represents the first $m$ terms, i.e. $f'(k, 0), f'(k, 1), \ldots, f'(k, m)$, that are common to every primitive recursive function such as $f'(x_1, x_2)$ that is formally represented by $F(x_1, x_2, x_3)$.

To prove this for any given $F(x_1, x_2, x_3)$, and given natural numbers $k$ and $m$, we can determine a suitable Turing machine that will actually *construct* the sequence $f'(k, 0), f'(k, 1), \ldots, f'(k, m)$ uniquely in order to verify the assertion.

Clearly we can choose a $k$ such that the number of terms in the above sequence, and so the Gödel-number of any Turing-computable proof, exceeds any given natural number $k_1$.

It follows that, in *omega$_2$*-**PA**, we cannot assume that if the Gödel-number of a proof of:

$\quad$ $(E_1x_3)F(x_1, x_2, x_3)$

is less than a given natural number $k_1$, then the Gödel-number of a proof of:

$\quad$ $(E_1x_3)F(x+1, x_2, x_3)$

will also not exceed $k_1$.

We cannot thus assume that:

$\quad$ *omega$_2$*-**PA** proves: $(E_1x_3)F(x_1, x_2, x_3) \Rightarrow (E_1x_3)F(x+1, x_2, x_3)$.

Hence we cannot assume that the argument of §**2.13** will establish strong representability for every recursive function $F(x_1, x_2, x_3)$, as is necessary for Gödel's reasoning to hold.

In fact, since the assumption that every recursive function can be *strongly* represented in any system of Arithmetic with the *omega-Constructivity* rule of inference **I**(*iii*) leads to inconsistency in the system, it follows that Gödel's conclusions do not logically follow from his reasoning in any interpretation of *omega$_2$*-**PA**.

We once more have the apparent conflict with the implicit Platonistic thesis that it is the axioms of a formal system that essentially determine the logical consequences in the interpretation.

**4.11. Does Gödel's reasoning hold even when his conclusions do not?**

The question arises: When can we infer Gödel's conclusions from his reasoning?

# Chapter 5. Gödel's argument for undecidable proposition

**Contents**

## 5.0. Introduction

We have argued in this paper that though Gödel's reasoning for the formal undecidability of true propositions is valid in all the above systems of Arithmetic, the inferences that can be validly drawn from it depend on the choice of the particular system.

We now highlight the conditional nature of Gödel's conclusions in a general Arithmetic **GA**.

We also argue that Gödel's undecidable formulas have the form of palimpsests that translate as ill-defined sentences in every interpretation

## 5.1. Every recursive function can be constructively represented in any Arithmetic

We note that, by a lemma of Hilbert and Bernays, if $f(x, y)$ is a recursive function, then it can be constructively represented in **GA** by the well-formed formula $F(x_1, x_2, x_3)$, as defined in §2.2 above.

Specifically, this means that, in any of the above systems, we can construct Turing machines to verify that, for every set of natural numbers $k, m, n$, if $f(k, m) = n$ then:

(i) **GA** proves: $F(\boldsymbol{k}, \boldsymbol{m}, \boldsymbol{n})$.

(ii) **GA** proves: $(Az)(\sim(z=\boldsymbol{n}) => \sim F(\boldsymbol{k}, \boldsymbol{m}, z))$.

## 5.2. Every representation of a recursive function has a unique Gödel-number

Since every recursive function, by definition, is built up from only a finite number of recursive functions of lower rank, we can further construct a Turing machine so that $F(x_1, x_2, x_3)$ can be expressed in only the primitive symbols of **GA**, and so uniquely coded and Gödel-numbered by some natural number $F\#$.

### 5.3. Gödel's recursive definition of provability

Such a process of coding allows us, following Gödel, to define a recursive arithmetic relation *prf*($x$, $y$), constructed out of 44 '*simpler*' recursive arithmetic functions and relations (*Gödel 1931 p182*), such that *prf*($x$, $y$) is constructively true if and only if $x$ is the Gödel-number of a finite proof sequence $X$ in **GA** for some well-formed formula $Y$ in **GA** whose Gödel-number is $y$.

In other words, for any given natural numbers $k$, $m$, we can construct a Turing machine *T*(*prf*($k$, $m$)) to write out the finite expression of **GA** whose Gödel-number is $k$, check if this is a valid proof sequence of **GA**, extract the last element of the sequence, and check if the Gödel-number of this element is $m$, returning the value *T*(*prf*($k$, $m$))=0 if *prf*($k$, $m$) is true and 1 if *prf*($k$, $m$) is false.

### 5.4. Gödel's Turing-computable *q*($x$, $y$)

We can now define a Turing-computable function *q*($x$, $y$) which is true if and only if $x$ is the Gödel-number of a well-formed formula *K*($z$) of **GA** with a single free variable $z$, and $y$ is the Gödel-number of a proof of *K*($x$) in **GA**.

In other words, given any numbers *K*#, *M*#, we can establish a Turing machine *T*(*K*#, *M*#) which will decode *K*# to yield the string *K*$ of **GA**, check if *K*$ is a well-formed formula in **GA** of the form *K*($z$), decode and check if the string *M*$ is a valid proof sequence in **GA**, and then check if the last well-formed formula of the sequence is *K*(*K*#).

If so, it establishes that *q*(*K*#, *M*#) is true; if not, it establishes that *q*(*K*#, *M*#) is false.

## 5.5. Gödel's constructive self-reference lemma

Hence the constructive self reference involved in Gödel's argument is:

$q(K\#, M\#)$ is true $\iff$ $M\$$ is a proof in **GA** of some well-formed formula $K(K\#)$ of **GA**.

## 5.6. Hilbert and Bernay's representation lemma

Now, by §**5.1**, since $q(x, y)$ is recursive, there is a well-formed formula $Q(x, y)$ of **GA** such that, for any natural numbers $K\#$, $M\#$:

(i) $q(K\#, M\#)$ is true $\Rightarrow$ $Q(K\#, M\#)$ is provable in **GA**, and

(ii) $q(K\#, M\#)$ is false $\Rightarrow$ $\sim Q(K\#, M\#)$ is provable in **GA**.

## 5.7. The Gödelian premise

We assume now that the well-formed formula $P\$$ given by $(Ay)(\sim Q(x, y))$, and the well-formed formula $G\$$ given by $(Ay)(\sim Q(P\#, y))$, are well-formed formulas of a simply consistent **GA**.

## 5.8. G$ is not provable

Assuming $G\$$ is provable in **GA**, let $M\$$ be the proof of $(Ay)(\sim Q(P\#, y))$ in **GA**.

We then have that $q(P\#, M\#)$ is true by definition.

Hence $Q(P\#, M\#)$ is provable in **GA**, and so $(Ey)Q(P\#, y)$ is provable in **GA**.

Hence $\sim(Ay)(\sim Q(P\#, y))$ is provable in **GA**, i.e. $\sim G\$$ is provable in **GA**

We conclude that (A*y*)(~*Q*(*P*#, *y*)) is not provable in **GA** if **GA** is simply consistent, i.e. *G*$ is not provable in **GA** if **GA** is simply consistent

**5.9. ~G$ is not provable**

Assuming ~*G*$ is provable in **GA**, it follows that *G*$ is not provable in **GA** if **GA** is simply consistent.

Hence *q*(*P*#, *y*) is false for all *y*.

It follows that ~*Q*(*P*#, *y*) is provable in **GA** for all *y*.

So, in every interpretation of **GA**, ~*Q*(*P*#, *y*) translates as a true proposition for all *y*.

However, since ~*G*$ is provable in **GA**, then ~(A*y*)(~*Q*(*P*#, *y*)) is provable in **GA**.

We then have the contradiction, in every interpretation of **GA**, that *Q*(*P*#, *M*#) translates as true for some *M*#, thus implying that **GA** is inconsistent.

We conclude that ~*G*$ is not provable in **GA** if **GA** is simply consistent.

**5.10. The essence of Gödel's reasoning**

What we have above is that, if *G*$ is a well-formed formula in **GA**, then:

    (*i*)    [(**GA** proves: *G*$) => (**GA** proves: ~*G*$)] => ~(**GA** proves: *G*$)

    (*ii*)    [(**GA** proves: ~*G*$) => (Every interpretation of **GA** is inconsistent)]
                                                                        => ~(**GA** proves: ~*G*$)

**5.11. Gödel's undecidable proposition GUS**

We consider now the well-formed formula *G*$ of **GA**, given by (A*y*)(~*Q*(*P*#, *y*)), under interpretation.

We argue below, in §**5.14**, that if $(Ay)(\sim Q(P\#, y))$ translates as a false sentence under any interpretation, then $Q(P\#, L\#)$ translates as a true assertion for some $L\$$.

However, we then have the contradiction that $(Ay)(\sim Q(P\#, y))$ translates as a true sentence under the interpretation.

We conclude that, if $G\$$ is a well-formed formula in **GA**, then $G\$$ translates into a true proposition in every interpretation of **GA**, even though both the well-formed formulas $G\$$ and $\sim G\$$ are not provable in **GA**, i.e. in every interpretation of **GA**, $G\$$ yields a Gödelian, formally undecidable but true, proposition **GUS**.

**5.12. Gödel's reasoning in PA**

Since, by §**2.13**, the well-formed formula $F(x_1, x_2, x_3)$ is well-defined in **PA** for any primitive recursive function $f(x, y)$, the above reasoning yields Gödel's proof of undecidability as a non-trivial assertion in **PA**.

**5.13. Gödel's reasoning does not establish his conclusions in most Arithmetics**

However, in all the other systems of Arithmetic defined above, Gödel's reasoning does not establish his conclusions significantly since the well-formed formula $F(x_1, x_2, x_3)$ is not obviously well-defined in these systems for every primitive recursive function $f(x, y)$.

In systems of Arithmetic with an *omega*-constructive rule, in particular, the assumption §**5.7** yields the meta-inconsistency:

(i)  $[(\mathbf{GA}\ \text{proves:}\ G\$) \Rightarrow (\mathbf{GA}\ \text{proves:}\ \sim G\$)] \Rightarrow \sim(\mathbf{GA}\ \text{proves:}\ G\$)$

(ii) $[\sim(\mathbf{GA}\ \text{proves:}\ G\$) \Rightarrow (\mathbf{GA}\ \text{proves:}\ G\$)] \Rightarrow (\mathbf{GA}\ \text{proves:}\ G\$)$

Ipso facto, both *G*$ and ~*G*$ are not well-defined, well-formed formulas in any systems of Arithmetic with an *omega*-constructive rule.[10]

They are thus trivially unprovable, and Gödel's formally undecidable proposition **GUS** trivially true, in every interpretation of a system of Arithmetic with an *omega*-constructive rule.

**5.14. The Palimpsest syndrome in standard PA : Overview of the argument**

We now argue that **GUS** essentially formalises the *Liar* sentence in standard **PA** by means of a "palimpsest".[11]

---

[10] We note that Russell's anomalous set **R**, defined by {*x* | '*x* is a member of **R**', '*x* is not a member of *x*'} as the set of all sets that are not members of themselves, can also be expressed in terms of the primitive symbols of an axiomatic set theory **S**. However, the assumption that the well-formed formulas '**R** is a member of **R**' and '**R** is not a member of **R**' are well-defined, well-formed formulas of **S** is inconsistent with the axioms of **S**.

[11] The intention here is to translate the formula '(A*y*)(~*Q*(*P*#, *y*))' from standard **PA** into English in order to argue the point that the truth of this assertion, under any interpretation, does not depend on the interpretation into which it is translated, but on the truth of a factual assertion, pertaining to standard **PA** specifically, that is embedded in the formula '(A*y*)(~*Q*(*P*#, *y*))' itself in the form of a "palimpsest".

For example, albeit at an explicit level, the statement "The English word 'tree' contains four letters of the English alphabet" must remain true no matter into which language - standard **PA**, Sanskrit, sign-language etc. - it is translated; it contains an assertion of fact that cannot be altered by the language into which it is translated.

If the destination-language contains the source-language, then it may indeed be able to determine the truth of factual assertions such as the above. If not, such assertions would be factually undecidable, at least until enough of the source-language is suitably adjoined.

The thesis here is that every interpretation of Gödel's formal system of standard **PA** in another language implicitly preserves the logical content of "named" entities; moreover, it subsumes standard **PA** and so contains the necessary wherewithal to determine the truth or falsity of factual assertions pertaining to the formal system of standard **PA**.

We argue that, by virtue of its roots in the finite construction, through Gödel-numbering, of the primitive recursive relation *q*(*P*#, *y*), every interpretation of Gödel's undecidable formula (A*y*)(**~Q**(*P*#, *y*)) can constructively be well-defined over the entire domain of any interpretation **M** in such a way that it is equivalent to the assertion:

> "For all *y*, the formula *Y*, whose Gödel-number is *y*, is not the interpretation of a finite proof sequence in standard **PA** for the formula (A*y*)(**~Q**(*P*#, *y*)), whose Gödel-number is *P*#, from the given set ***Axm*** of recursively defined axioms, and the given set ***Inf*** of recursively defined rules of inference, of standard **PA**."

If standard **PA** is assumed consistent, this means either that the sentence (A*y*)(**~Q**(*P*#\*, *y*)) (*where P#\* is the interpretation in M of the numeral P#*) can be interpreted as true in every interpretation **M** of standard **PA** (*since the set **Axm** of axioms, and set **Inf** of rules of inference, of Gödel's formal system of standard **PA**, are defined recursively*), or the expression (A*y*)(**~Q**(*P*#\*, *y*)) treated as an ill-defined sentence in every **M** under interpretation.

In other words, even if (A*y*)(**~Q**(*P*#, *y*)) is a well-defined formula in standard **PA**, there can be no consistent, non-standard model **M** of standard **PA** in which its interpretation (A*y*)(**~Q**(*P*#\*, *y*)) is uniquely false.

> (The latter would imply the contradiction that we can write out a finite proof sequence in standard **PA** for the unprovable formula (A*y*)(**~Q**(*P*#, *y*)) using only the given set ***Axm*** of axioms and the given set ***Inf*** of rules of inference of standard **PA**.)

By virtue of Gödel's completeness theorem for first order predicate calculus, it then follows that even if, like Russell's paradoxical set, the formula (A*y*)(**~Q**(*P*#, *y*)) appears to be a valid construction in standard **PA**, it must translate as an ill-defined sentence under every interpretation.

In other words, there is an essential element in the definition of the formula (A*y*)(~*Q*(*P*#, *y*)) such that, as in the case of the ***Liar*** sentence, we cannot constructively ascribe, or even non-constructively assume, either truth or falsity to/for the sentence (A*y*)(~*Q*(*P*#*, *y*)) in an interpretation without inconsistency.

It is a moot point whether we choose to treat (A*y*)(~*Q*(*P*#*, *y*)) as an ill-defined "sentence" in every interpretation, or to treat standard **PA** as an inconsistent system.

### 5.15. The Extended Representation Lemma

To give the argument in detail, we note the following meta-properties of recursive arithmetical relations *f*(*y*) where, if **M** is any interpretation of standard **PA** then, for all natural numbers *n*, and their representations ***n*** and interpretations ***n***\* in standard **PA** and **M** respectively:

    *f*(*n*) is true in the standard interpretation

        => ***F***(***n***) is provable in standard **PA**

            => ***F***(***n***\*) is true in the interpretation **M**

    *f*(*n*) is false in the standard interpretation

        => ~***F***(***n***) is provable in standard **PA**

            => ***F***(***n***\*) is false in the interpretation **M**

From this we conclude that, for every natural number *n*:

    *f*(*n*) is true in the standard interpretation  <=>  ***F***(***n***\*) is true in the interpretation **M**

If we extend the definition of *f*(*x*) in the standard interpretation **SI** so that *f*(*x*) is false if *x* is not a natural number, then we have the Extended Representation meta-Lemma for recursive arithmetical functions *f*(*x*):

*f*(*x*) is true in the standard interpretation  <=>  *F*(*x*) is true in the interpretation **M**

## 5.16. The Palimpsest syndrome in standard PA : The proof

We note that Gödel's undecidable formula ~(A*y*)(~*Q*(*P*#, *y*)) now translates as a well-defined sentence over the entire domain of any interpretation **M**.

The significance of the Extended Representation meta-Lemma is that, by virtue of Goedel's original definition of the primitive recursive relation *q*(*P*#, *y*) - involving only finite Gödel-numbers of finite sets of well-formed, well-defined finite formulas and finite proof sequences of standard **PA** - its representation *Q*(*P*#, *y*) in standard **PA** is well-defined only over the denumerable domain of the numerals *n* in standard **PA**, corresponding to the denumerable natural numbers *n*.[12]

---

[12] The central issue in these papers actually concerns the relationship of the formally represented elements to the elements that they seek to represent. Thus, if standard **PA** is to serve as the formal representation of intuitive arithmetic **IA**, then the questions that must be addressed are whether every element of intuitive arithmetic **IA** is represented in standard **PA**, and whether there are elements of standard **PA** that are not reflected in intuitive arithmetic **IA**.

If we assume that our intuitive arithmetic **IA** is what we define as the standard interpretation **SI** of standard **PA**, then the question is whether every element of the standard interpretation **SI** is represented in standard **PA**, and whether there are elements of standard **PA** that are not reflected in the standard interpretation **SI**.

In the former case, standard **PA** must be judged as inadequate to our intent; in the latter, as containing elements extraneous to our intent.

Applying the same measure to the representation of, say, a function *f*(*x*) of the standard interpretation **SI** by some formula *F*(*x*, *y*) of standard **PA**, it follows that if *f*(*x*) is an intuitively well-defined function, then *F*(*x*, *y*) too must be a well-defined representation,

It follows that, in any interpretation **M** of standard **PA**, the sentence ~(A*y*)(**~*Q*(*P*#\*, *y*)) - which is equivalent to the assertion (E*y*)*Q*(*P*#\*, *y*) - is well-defined if and only if the relation *Q*(*P*#\*, *a*\*) is satisfied for some value *a*\* in the domain of **M** that is an interpretation of some numeral *a* of standard **PA**. However, we cannot conclude from this that *Q*(*P*#\*, *a*\*) is well-defined - and so decidable - for every value *a*\* in the domain of **M** that is not an interpretation of some numeral *a* of standard **PA**.

This lacuna is now filled by the Extended Representation meta-Lemma, which ensures that the range over which the relation *Q*(*P*#\*, *a*\*) is well-defined covers the entire domain of **M**, and that *Q*(*P*#\*, *a*\*) is not satisfied by any value *a*\* in the domain of **M** that is not an interpretation of some numeral *a* of standard **PA**.

This follows immediately since (A*y*)(**~*Q*(*P*#\*, *y*)) in **M** is now equivalent to:

> "For all *y*, the formula *Y*, whose Gödel-number is *y*, is not the interpretation in **M** of a finite proof sequence in standard **PA** for the formula (A*y*)(**~*Q*(*P*#, *y*)), whose Gödel-number is *P*#, from the given set ***Axm*** of recursively defined axioms, and the given set ***Inf*** of recursively defined rules of inference, of standard **PA**."

---

i.e. *F*(*x*, *y*) cannot be considered to represent *f*(*x*) in some respects only, but must be taken to represent it entirely.

We argue that Gödel's interpretation of undecidability actually violates this basic integrity. It implies that the truth of some propositions cannot be reflected in the formalisation, whilst the undecidability of some formulas cannot be reflected in the assertions that they represent.

We argue that this interpretation is subjective, and that a more objective interpretation implies that the undecidability of some formulas in standard **PA** reflects assertions in the standard interpretation **SI** that cannot be ascribed any truth-value, and are thus ill-defined assertions similar to the *Liar* sentence, and not true sentences as concluded by Gödel.

In other words, the sentence ~(A*y*)(**~Q**(***P*#***, *y*)) is well-defined in **M**, and ***Q*(*P*#***, ***a****) is necessarily false for all values ***a**** in the domain of **M** that are not interpretations of some numeral ***a*** of standard **PA**.

> (The argument here is that the relation ***Q***(*x*, *y*) can always be well-defined in every interpretation **M** of standard **PA** in such a way that it mirrors the "palimpsest" that is embedded in the definition of the relation ***q***(*x*, *y*) of the standard interpretation **SI**.)

We thus argue that there is no interpretation in which the formula ~(A*y*)(**~Q**(***P*#**, *y*)) of standard **PA** uniquely translates as a well-defined sentence ~(A*y*)(**~Q**(***P*#***, *y*)) that can be true for some value ***a**** in the domain of **M** that is not an interpretations of some numeral ***a*** of standard **PA**.

Hence the truth of the sentence ~(A*y*)(**~Q**(***P*#***, *y*)) in any interpretation **M** of standard **PA** can always be defined so as to imply the truth of ***Q*(*P*#***, ***n****) for some ***n**** in **M**, where ***n**** is the interpretation of the numeral ***n*** in **M**.

This, of course, implies the truth of ***q***(***P*#**, *n*) for some natural number *n* in the standard interpretation **SI** by the Extended Representation meta-Lemma. This in turn, by Gödel's Self-reference Lemma for ***q***(***P*#**, *y*), implies that *n* is a proof of the sentence (A*y*)(**~Q**(***P*#**, *y*)) in standard **PA**, contradicting the unprovability of the formula (A*y*)(**~Q**(***P*#**, *y*)) in standard **PA**.

Hence the sentence ~(A*y*)(**~Q**(***P*#***, *y*)) can be falsified in every interpretation **M** of standard **PA**, so that the sentence (A*y*)(**~Q**(***P*#***, *y*)) is true in every interpretation **M** of standard **PA**.

However, by Gödel's completeness theorem for first order predicate logic, this again implies that the formula (A*y*)(**~Q**(***P*#**, *y*)) is provable in standard **PA**.

## 5.17. Conclusion

The contradiction establishes that, if standard **PA** is consistent, then Gödel's undecidable formulas, (A*y*)(~*Q*(*P*#, *y*)) and ~(A*y*)(~*Q*(*P*#, *y*)), translate as ill-defined sentences under interpretation[13].

# Chapter 6. Conclusions

**Contents**

---

[13] The cause of such ill-definedness apparently lies in the formal definition of *Q*(*x*, *y*) as the representation of *q*(*x*, *y*) in standard **PA** through use of the Gödel Beta function. By definition, *Q*(*x*, *y*) thus represents, in standard **PA**, an infinity of "recursive" functions of every interpretation (*cf. §2.9*).

Hence the **PA** provability of *Q*(*k*, *m*) may be intuitively viewed as corresponding, in any interpretation **M**, to the intuitive assertion that, for any interpretations *k*\* and *m*\* in **M** of the numerals *k* and *m*, we can always *construct* some (*non-unique*) pair of values *u*\*, *v*\* such that we can define a "Beta-function" *ß*\*(*u*\*, *v*\*, *i*) in **M** that represents the first *m*\* terms, i.e. *q'*(*k*\*, 0\*), *q'*(*k*\*, 1\*), ... , *q'*(*k*\*, *m*\*), which are common to an infinity of "primitive recursive relations" of **M**, such as *q'*(*x*, *y*), that are formally represented in **PA** by *Q*(*x*, *y*)..

Clearly, in this view, the provability of the formula *Q*(*x*, *y*) for variable *x* and *y* does not interpret similarly into a unique, well-defined proposition that can be asserted as true in **M**. It follows that the provability of ~*Q*(*x*, *y*) in **PA** interprets as a true proposition in **M** by virtue of its ill-definedness (*i.e. although we may not intuitively assert an ill-defined proposition as true in **M**, we may yet assert that it is not a true proposition in **M***).

### 6.1. The non-constructive nature of Gödel's first order Arithmetic

We have argued in §**1** and §**2** that the non-constructive nature of Gödel's system of standard first order Arithmetic **PA** admits *omega*-provable well-formed formulas that may translate into propositions that refer to ill-defined totalities (*interpreted Platonistically*) such as those involved in the logical antinomies.

### 6.2. *Simply* constructive first order systems

In §**3**, we defined various systems of first order Arithmetic in which well-formed formulas that may correspond to ill-defined predicates under interpretation cannot be strongly represented if the systems are simply consistent. We argued, accordingly, that Gödel's reasoning does not hold in most simply consistent systems of Arithmetic.

We defined a *simply* constructive system, *omega*-**PA**, of first order Peano's Arithmetic **PA**, and argued that every model of **PA** is a model of a *simply* constructive *omega*-**PA**. We argued further that *omega*-**PA**, unlike **PA**, admits Cantor's *ordinal*-Arithmetic **CA** as a model.

### 6.3. *Strongly* constructive first order systems

In §**4**, we extended the scope of the above argument by defining *stronger* constructive systems, *omega$_1$*-**PA** and *omega$_2$*-**PA**, of first order Peano's Arithmetic **PA**.

We argued that, even here, well-formed formulas that may correspond to ill-defined predicates under interpretation cannot be strongly represented if the system is simply consistent, and so again Gödel's reasoning does not hold.

However, we have also argued that the axioms, interpretations and models of **PA** and the axioms, interpretations and models of a *strongly* constructive *omega$_2$*-**PA** are identical.

**6.4. The roots of Gödel's Platonism**

We thus suggest that Gödel's Platonism may not have been the result of a faith-only belief in an abstract world of absolute mathematical ideals.

If the arguments of this paper are substantive, then Gödel's Platonism may have been an intuitive reflection of the non-constructive, and implicitly Platonistic, nature of the first order predicate calculus chosen by him for defining his formal system of Arithmetic **PA**.

**6.5. The conflict within model theory**

We have also highlighted the conflict with classical model theory, and argued that the consequences within a model are not Platonistically determined absolutely by the interpretation of only the axioms of a formal system in the model, but are consequences of the rules of inference that we select for the system.

**6.6. Consistency**

We conclude with the tentative remarks that, from the above, it can be reasonably argued that consistency may not be an inherent feature of an axiomatic system.

It may, instead, be viewed as a feature of specification that we design into the definition of a system through our choice of appropriate rules of inference in order to specify those elements of the system, and its interpretations, that reflect what we intend the system to communicate faithfully.

**6.7. Gödel's undecidable proposition is an ill-defined sentence**

We argue in §5 that Gödel's undecidable formula **GUS**, if assumed well-defined under interpretation, translates as a true sentence in every interpretation of standard **PA**. It follows that there can be no consistent, non-standard, interpretation in which **GUS** is false.

Since, by virtue of Gödel's completeness theorem for a first order predicate calculus, the above argument contradicts the essential unprovability of **GUS** in standard **PA**, we argue further that, if standard **PA** is assumed simply consistent, then Gödel's undecidable formula **GUS** translates as an ill-defined sentence that mirrors the *Liar* sentence in every interpretation of standard **PA**.

## References

**Burbanks, Andrew. (*e-seminar*) *Fast-Growing Functions and Unprovable Theorems*.**
<Home page : http://www.maths.bris.ac.uk/~maadb/>
<Seminar : http://www.maths.bris.ac.uk/~maadb/research/seminars/online/fgfut/index.html>
<Ord'ls : http://www.maths.bris.ac.uk/~maadb/research/seminars/online/fgfut/fgfut09.html>

**Davis, M (ed.). 1965. The Undecidable. Raven Press, New York.**
<Home page : http://www.cs.nyu.edu/cs/faculty/davism/>

**Franzén, Torkel. *Provability and Truth*.**
< Home page : http://www.sm.luth.se/~torkel/>
<See also remarks : http://www.sm.luth.se/~torkel/eget/godel/system.html>
<See also remarks : http://www.sm.luth.se/~torkel/eget/thesis/chapter5.html>

**Friedman, Harvey. 1997. *Unprovable theorems in discrete mathematics*.**
<Unprovable : http://www.math.berkeley.edu/~ribet/Colloquium/friedman.html>

**Gödel, Kurt. 1931. *On formally undecidable propositions of Principia Mathematica and related systems I*. Also in M. Davis (ed.). 1965. The Undecidable. Raven Press, New York.**
<1931 Paper : http://www.ddc.net/ygg/etext/godel/godel3.htm>
<Formal system : http://www.ddc.net/ygg/etext/godel/godel3.htm - 15>
<Reducibility : http://www.ddc.net/ygg/etext/godel/godel3.htm - AXIOM-IV>
<Comprehension : http://www.ddc.net/ygg/etext/godel/godel3.htm - AXIOM-IV>
<See remarks on p190/191 : http://www.ddc.net/ygg/etext/godel/godel3.htm - 16>
<Recursive definitions : http://www.ddc.net/ygg/etext/godel/godel3.htm - DEF1>
< Beta-function : http://www.ddc.net/ygg/etext/godel/godel3.htm - LEMMA1>
< See remarks on p189 : http://www.ddc.net/ygg/etext/godel/godel3.htm - 16>
<Eqn. 8.1 : http://www.ddc.net/ygg/etext/godel/godel3.htm - EQ8.1>
<Axiom IV : http://www.ddc.net/ygg/etext/godel/godel3.htm - AXIOM-IV>
<Footnote 48a : http://www.ddc.net/ygg/etext/godel/godel3.htm - 48a>
<cf. Recursive definition #17 : http://www.ddc.net/ygg/etext/godel/godel3.htm - DEF17>
<cf. Proposition V : http://www.ddc.net/ygg/etext/godel/godel3.htm#PROP-V>

**Gödel, Kurt. 1934.** *On undecidable propositions of formal mathematical systems*. In M. Davis (ed.). 1965. The Undecidable. Raven Press, New York.

**Lucas, J. R.** *The Gödelian argument*.
    <Home page : http://users.ox.ac.uk/~jrlucas/>
    <Gödelian argument : http://www.leaderu.com/truth/2truth08.html>

**Mendelson, Elliott. 1964. Introduction to Mathematical Logic. Van Norstrand, Princeton.**
    <Home page : http://sard.math.qc.edu/Web/Faculty/mendelso.htm>

**Parikh, Rohit. 1973.** *On the length of proofs*. In Trans. AMS, 177, 29-36
    <Home page : http://www.sci.brooklyn.cuny.edu/cis/parikh/>

**Parsons, Charles D. 1995.** *Platonism and mathematical intuition in Kurt Gödel's thought*. **In The Bulletin of Symbolic Logic, Volume 1, Issue 1, March 1995.**
    <1995 Article : http://www.math.ucla.edu/~asl/bsl/0101-toc.htm>

**Peregrin, Jaroslav. 1997.** *Language and its Models: Is Model Theory a Theory of Semantics?*
    < cf. also Models : http://www.hf.uio.no/filosofi/njpl/vol2no1/models/models-html.html>
    <cf. also Quantifier : http://www.hf.uio.no/filosofi/njpl/vol2no1/models/node8.html>
    <cf. also Basis : http://www.hf.uio.no/filosofi/njpl/vol2no1/models/node5.html>
    < cf. also Consequence : http://www.hf.uio.no/filosofi/njpl/vol2no1/models/node11.html>
    < cf. also Truth : http://www.hf.uio.no/filosofi/njpl/vol2no1/models/node3.html>
    < cf. also Language : http://www.hf.uio.no/filosofi/njpl/vol2no1/models/node2.html>

**Podnieks, Karlis. 2001. Around Goedel's Theorem.**
    <Home page : http://www.ltn.lv/~podnieks/index.html>
    <Around Goedel's Theorem : http://www.ltn.lv/~podnieks/gt.html>
    <Goedel's Incompleteness Theorem : http://www.ltn.lv/~podnieks/gt5.html - BM5_3>
    <Standard first order Arithmetic **PA** : http://www.ltn.lv/~podnieks/gt3.html - BM3>
    <Formula-number : http://www.ltn.lv/~podnieks/gt5.html - BM5_2>
    <Decode : http://www.ltn.lv/~podnieks/gt5.html - BM5_2>
    <Representation theorem : http://www.ltn.lv/~podnieks/gt3.html - BM3_3>
    <Platonism : http://www.ltn.lv/~podnieks/gt1.html>
    <Self-reference lemma : http://www.ltn.lv/~podnieks/gt5.html - BM5_2>

**Rosmaita, Brian J.** *Minds, machines, and metamathematics: constraints on the mathematical objection to mechanism*.
    <Article : http://www.personal.kent.edu/~brosmait/apa96.html>

## Web resources

*Free On-line Dictionary of Computing*.
    <Dictionary : http://www.nightflight.com/foldoc/>
    <cf. constructive : http://lgxserve.ciseca.uniba.it/lei/foldop/foldoc.cgi?query=constructive>
    <Simply consistent : http://www.swif.uniba.it/lei/foldop/foldoc.cgi?simple+consistency>
    <Infinite : http://www.swif.uniba.it/lei/foldop/foldoc.cgi?query=infinite&action=Search>
    <Generalisation : http://www.swif.uniba.it/lei/foldop/foldoc.cgi?generalization>
    <Converse : http://lgxserve.ciseca.uniba.it/lei/foldop/foldoc.cgi?query=converse>

<Turing machine : http://www.swif.uniba.it/lei/foldop/foldoc.cgi?query=turing>
<Artificial Int. : http://www.swif.uniba.it/lei/foldop/foldoc.cgi?query=AI&action=Search>
<Range : http://www.nightflight.com/foldoc-bin/foldoc.cgi?query=range>
<Cmprhn. axiom : http://ase.isu.edu/ase01_07/ase01_07/bookcase/ref_sh/foldoc/99/8.htm>
<ZF set theory : http://ase.isu.edu/ase01_07/ase01_07/bookcase/ref_sh/foldoc/35/122.htm>
<Mapping : http://www.swif.uniba.it/lei/foldop/foldoc.cgi?query=mapping&action=Search>

*Free On-line Dictionary of Philosophy.*
<Dictionary : http://www.swif.uniba.it/lei/foldop/>

*Glossary of First-Order Logic.*
<Glossary : http://www.earlham.edu/~peters/courses/logsys/glossary.htm>
<First order : http://www.earlham.edu/~peters/courses/logsys/glossary.htm - fot>
<Predicate calculus : http://www.earlham.edu/~peters/courses/logsys/glossary.htm - pl>
<Arithmetic : http://www.earlham.edu/~peters/courses/logsys/glossary.htm - arithmetic>
<Omega : http://www.earlham.edu/~peters/courses/logsys/glossary.htm - omegacom>
<Formula : http://www.earlham.edu/~peters/courses/logsys/glossary.htm - formalsystem>
<Domain : http://www.earlham.edu/~peters/courses/logsys/glossary.htm - domain>
<Set : http://www.earlham.edu/~peters/courses/logsys/glossary.htm - set>
<Natural numbers : http://www.earlham.edu/~peters/courses/logsys/glossary.htm - naturals>
<Numeral : http://www.earlham.edu/~peters/courses/logsys/glossary.htm - numeral>
<Compound : http://www.earlham.edu/~peters/courses/logsys/glossary.htm - compoundprop>
<cf. Meta... : http://www.earlham.edu/~peters/courses/logsys/glossary.htm - metalanguage>
<Intrprttn. : http://www.earlham.edu/~peters/courses/logsys/glossary.htm - interpretation>
<Equivalent : http://www.earlham.edu/~peters/courses/logsys/glossary.htm - equivalence>
<Omega-comp. : http://www.earlham.edu/~peters/courses/logsys/glossary.htm - omegacom>
<Free variable : http://www.earlham.edu/~peters/courses/logsys/glossary.htm - freevar>
<Consequence : http://www.earlham.edu/~peters/courses/logsys/glossary.htm - syncon>
<Proof sequence : http://www.earlham.edu/~peters/courses/logsys/glossary.htm - proof>
<Denum. : http://www.earlham.edu/~peters/courses/logsys/glossary.htm - denumerableset>
<Rcrsv. : http://www.earlham.edu/~peters/courses/logsys/glossary.htm - recursivefunction>
<Frml sys. : http://www.earlham.edu/~peters/courses/logsys/glossary.htm - formallanguage>
<Cptbl. : http://www.earlham.edu/~peters/courses/logsys/glossary.htm - computablefunction>
<Standard model : http://www.earlham.edu/~peters/courses/logsys/glossary.htm - model>
<Set theory : http://www.earlham.edu/~peters/courses/logsys/glossary.htm - settheory>
<Closure : http://www.earlham.edu/~peters/courses/logsys/glossary.htm - closure>
<Model : http://www.earlham.edu/~peters/courses/logsys/glossary.htm - model>
<Std. set theory : http://www.earlham.edu/~peters/courses/logsys/glossary.htm - settheory>
<Epimenides : http://www.earlham.edu/~peters/courses/logsys/glossary.htm - enumset>
<Russell : http://www.earlham.edu/~peters/courses/logsys/glossary.htm - russellparadox>
<Gödel-number : http://www.earlham.edu/~peters/courses/logsys/glossary.htm#gnumbering>
<Well-formed formula : http://www.earlham.edu/~peters/courses/logsys/glossary.htm - wff>
<Rprsntn. : http://www.earlham.edu/~peters/courses/logsys/glossary.htm - repfunction>
<Proposition : http://www.earlham.edu/~peters/courses/logsys/glossary.htm - proposition>
<Predicate : http://www.earlham.edu/~peters/courses/logsys/glossary.htm - powerset>
<Rules of inf. : http://www.earlham.edu/~peters/courses/logsys/glossary.htm - rulesinf>
<Total fn. : http://www.earlham.edu/~peters/courses/logsys/glossary.htm - totalfunction>
<Satisfaction : http://www.earlham.edu/~peters/courses/logsys/glossary.htm - satisfaction>

<Axioms : http://www.earlham.edu/~peters/courses/logsys/glossary.htm - axioms>
<Function : http://www.earlham.edu/~peters/courses/logsys/glossary.htm - function>
<Dcdbl : http://www.earlham.edu/~peters/courses/logsys/glossary.htm - decidablesystem>
<Transfinite : http://www.earlham.edu/~peters/courses/logsys/glossary.htm - totalfunction>
<Exstnc. pr. : http://www.earlham.edu/~peters/courses/logsys/glossary.htm - existenceproof>
<Cnst. pr. : http://www.earlham.edu/~peters/courses/logsys/glossary.htm - constructiveproof>
<Implication : http://www.earlham.edu/~peters/courses/logsys/glossary.htm - implication>
<Completeness : http://www.earlham.edu/~peters/courses/logsys/glossary.htm - completeness>

*Internet Encyclopedia of Philosophy*
<Antinomies : http://www.utm.edu/research/iep/p/par-log.htm>

*Larry Hauser's Mostly Modern Philosophical Glossary.*
<Glossary : http://members.aol.com/lshauser2/lexicon.html>
<True : http://members.aol.com/lshauser2/lexicon.html - truth>

*MacTutor History of Mathematics archive*
<Home page : http://www-groups.dcs.st-and.ac.uk/~history/>
<Kurt Gödel : http://www-groups.dcs.st-and.ac.uk/~history/Mathematicians/Godel.html >

*MathPages*
<MathPages : http://www.mathpages.com/home/>

*Stanford Encyclopedia of Philosophy.*
<Encyclopedia : http://setis.library.usyd.edu.au/stanford/archives/win1997/contents.html>
<Vagueness : http://setis.library.usyd.edu.au/stanford/archives/win1997/entries/vagueness/>
<cf. also Non-constructive : http://plato.stanford.edu/entries/mathematics-constructive/ - 2>
<Intuitionistically : http://plato.stanford.edu/entries/logic-intuitionistic/>
<Principia : http://plato.stanford.edu/entries/principia-mathematica/>
<cf. also Constructive : http://plato.stanford.edu/entries/mathematics-constructive/ - 1>
<cf. also Platonism : http://plato.stanford.edu/entries/mathphil-indis/ - 1>

*Web Dictionary of Cybernetics and Systems*
<Dictionary : http://pespmc1.vub.ac.be/ASC/indexASC.html>

*xrefer*
<Home page : http://www.xrefer.com/>
<Axiom of reducibility : http://www.xrefer.com/entry/553361>
<cf. Formalism : http://www.xrefer.com/entry/552091>

*Author's e-mail*: anandb@vsnl.com

(*Updated : Monday 28$^{th}$ January 2002 8:39:04 PM by *re@alixcomsi.com*)